\documentclass{article} 

\pdfoutput=1  %  for arxiv

\usepackage{enumerate,wrapfig,placeins}

\usepackage{amsmath}

\usepackage{caption}
\newtheorem{theorem}{Theorem}[section]
\newtheorem{lemma}[theorem]{Lemma}
\newtheorem{proposition}[theorem]{Proposition}

\newlength{\noteWidth}
\long\def\notes#1{\ifinner
	{\tiny #1}
\else
\marginpar{\parbox[t]{\noteWidth}{\raggedright\tiny #1}}
\fi}
\def\notes#1{}

\usepackage{textgreek,upgreek,bm}

\usepackage{titlesec}

\usepackage[usenames, dvipsnames]{color}

\usepackage[pdfa,hidelinks]{hyperref}

\usepackage{optidef}
\usepackage{graphics} % for pdf, bitmapped graphics files
\usepackage{subfigure}

\graphicspath{{../figures/}}
\usepackage{amsmath} % assumes amsmath package installed
\usepackage{amssymb}  % assumes amsmath package installed
\usepackage{balance}

\usepackage{verbatim}

\usepackage{geometry} 

\def\urls#1{{\footnotesize\url{#1}}}

\def\mindex#1{\index{#1}}

\DeclareFontFamily{U}{mathx}{\hyphenchar\font45}
\DeclareFontShape{U}{mathx}{m}{n}{<-> mathx10}{}
\DeclareSymbolFont{mathx}{U}{mathx}{m}{n}
\DeclareMathAccent{\widebar}{0}{mathx}{"73}

\makeatletter
\newcommand\gobblepars{%
    \@ifnextchar\par%
        {\expandafter\gobblepars\@gobble}%
{}}
\makeatother

\def\whamrm#1{\smallbreak\pagebreak[3]%
	\noindent\text{\rm#1}\ \ \gobblepars}
	
\def\whamit#1{\smallbreak\pagebreak[3]%
	\noindent\textit{#1}\ \ \gobblepars}

\def\wham#1{\smallbreak\pagebreak[3]%
	\noindent\textup{\textbf{#1}}\ \ \gobblepars}

\newcommand{\qsaprobe}{{\scalebox{1.1}{$\upxi$}}}  %\upzeta \textphi
\newcommand{\bfqsaprobe}{{\scalebox{1.1}{$\bm{\upxi}$}}}  

\def\odestate{\upvartheta}

\def\fee{\upphi}
\def\feex{\widetilde{\fee}}

\def\uQ{\underline{Q}}

\def\disc{\gamma}

\newcommand{\bbblot}{\raise1pt\hbox{\vrule height .4ex width .4ex depth .05ex}}
\long\def\defbox#1{\framebox[.9\hsize][c]{\parbox{.85\hsize}{%
\parindent=0pt
\baselineskip=12pt plus .1pt      % STYLE 
\parskip=6pt plus 1.5pt minus 1pt % CHANGES
 #1}}}

\long\def\beginbox#1\endbox{\subsection*{}%
\hbox{\hspace{.05\hsize}\defbox{\medskip#1\bigskip}}%
\subsection*{}}

\def\endbox{}

 \def\archival#1{} %Notes I'd like to save

\def\FRAC#1#2#3{\genfrac{}{}{}{#1}{#2}{#3}}

\def\ddt{{\mathchoice{\FRAC{1}{d}{dt}}%
{\FRAC{1}{d}{dt}}%
{\FRAC{3}{d}{dt}}%
{\FRAC{3}{d}{dt}}}}

\def\ddtp{{\mathchoice{\FRAC{1}{d^{\hbox to 2pt{\rm\tiny +\hss}}}{dt}}%
{\FRAC{1}{d^{\hbox to 2pt{\rm\tiny +\hss}}}{dt}}%
{\FRAC{3}{d^{\hbox to 2pt{\rm\tiny +\hss}}}{dt}}%
{\FRAC{3}{d^{\hbox to 2pt{\rm\tiny +\hss}}}{dt}}}}

\def\ddyp{{\mathchoice{\FRAC{1}{d^{\hbox to 2pt{\rm\tiny +\hss}}}{dy}}%
{\FRAC{1}{d^{\hbox to 2pt{\rm\tiny +\hss}}}{dy}}%
{\FRAC{3}{d^{\hbox to 2pt{\rm\tiny +\hss}}}{dy}}%
{\FRAC{3}{d^{\hbox to 2pt{\rm\tiny +\hss}}}{dy}}}}

\def\half{{\mathchoice{\FRAC{1}{1}{2}}%
{\FRAC{1}{1}{2}}%
{\FRAC{3}{1}{2}}%
{\FRAC{3}{1}{2}}}}

\def\darrow{\buildrel{\rm d}\over\longrightarrow}

\def\eqdist{\buildrel{\rm dist}\over =}

\def\argmin{\mathop{\rm arg{\,}min}}

\def\state{{\sf X}}

\def\ustate{{\sf U}} 
\def\ystate{{\sf Y}}

\def\bx{{{\cal B}(\state)}}

\def\ystate{{\sf Y}}

\def\by{\clB(\ystate)}

\def\bx{{{\cal B}(\state)}}

\def\bfmath#1{{\mathchoice{\mbox{\boldmath$#1$}}%
{\mbox{\boldmath$#1$}}%
{\mbox{\boldmath$\scriptstyle#1$}}%
{\mbox{\boldmath$\scriptscriptstyle#1$}}}}

\def\bfPhi{\bfmath{\Phi}}

\def\bfmI{\bfmath{I}}

\def\bfmU{\bfmath{U}}

\def\bfmX{\bfmath{X}}

\def\bfmY{\bfmath{Y}}

\def\bfmhhaY{\bfmath{\hhaY}} %\widehat{\widehat{Y}}}}
\def\bfmhhaY{\hbox to 0pt{$\widehat{\bfmY}$\hss}\widehat{\phantom{\raise 1.25pt\hbox{$\bfmY$}}}}

\def\bfmW{\bfmath{W}}

\def\hatheta{{\hat\theta}}

\def\tiltheta{{\tilde \theta}}

\def\clB{{\cal B}}
\def\clC{{\cal C}}
\def\clE{{\cal E}}

\def\clH{{\cal H}}
\def\clJ{{\cal J}}

\def\clS{{\cal S}}
\def\clX{{\cal X}}

\def\clZ{{\cal Z}}

\def\clC{{\cal C}}

\def\eqdef{\mathbin{:=}}

\def\Prob{{\sf P}}

\def\Expect{{\sf E}}

\def\lgmath#1{{\mathchoice{\mbox{\large #1}}%
{\mbox{\large #1}}%
{\mbox{\tiny #1}}%
{\mbox{\tiny #1}}}}

\def\Zero{{\mathchoice{\lgmath{\sf 0}}%
{\mbox{\sf 0}}%
{\mbox{\tiny \sf 0}}%
{\mbox{\tiny \sf 0}}}}

\def\ind{\bbbone}
  
 \def\epsy{\varepsilon}

\def\varble{\,\cdot\,}

\def\formtmp#1#2{{\vskip12pt\noindent\fboxsep=0pt\colorbox{#1}{\vbox{\vskip3pt\hbox to \textwidth{\hskip3pt\vbox{\raggedright\noindent\textbf{#2\vphantom{Qy}}}\hfill}\vspace*{3pt}}}\par\vskip2pt%
\noindent\kern0pt}}

\def\barf{{\widebar{f}}}
\def\barg{{\widebar{g}}}

\def\barA{{\bar{A}}}

\def\barJ{{\bar{J}}}

\def\bartheta{{\overline{\theta}}}

{\end{list}}

\def\ass(#1:#2){(#1\ref{#1:#2})}

\def\ritem#1{
\item[{\sf \ass(\current_model:#1)}]
}

\newenvironment{recall-ass}[1]{% 
\begin{description}
\def\current_model{#1}}{
\end{description}
}

 \newcommand{\blot}{\vrule height 1.1ex width .9ex depth -.1ex }
\def\qedb{\ifmmode\blot\else{\vspace{-.2cm}\unskip\nobreak\hfil
\penalty50\hskip1em\null\nobreak\hfil\blot
\parfillskip=0pt\finalhyphendemerits=0\endgraf}\fi}

  \makeatletter
\DeclareRobustCommand{\sqcdot}{\mathbin{\mathpalette\morphic@sqcdot\relax}}
\newcommand{\morphic@sqcdot}[2]{%
  \sbox\z@{$\m@th#1\centerdot$}%
  \ht\z@=.33333\ht\z@
  \vcenter{\box\z@}%
}
\makeatother

\newcommand{\qedsymbol}{\hbox{\tiny$\blacksquare$}}  %

\def\qed{\ifmmode\qedsymbol\else{\unskip\nobreak\hfil
\penalty50\hskip1em\null\nobreak\hfil\qedsymbol
\parfillskip=0pt\finalhyphendemerits=0\endgraf}\fi}

\newcounter{rmnum}

\newcounter{anum}

\newcommand{\field}[1]{\mathbb{#1}} 

\def\posRe{\field{R}_+}
\def\Re{\field{R}}

\def\Prob{{\sf P}}

\def\Expect{{\sf E}}

\def\transpose{{\intercal}}

\def\argmin{\mathop{\rm arg\, min}}
\def\ind{\hbox{\large \bf 1}}

\def\epsy{\varepsilon}
\def\varble{\,\cdot\,}

\def\haY{\widehat{Y}}

\def\hhaY{\hbox to 0pt{$\haY$\hss}\widehat{\phantom{\raise 1.25pt\hbox{Y}}}}

\def\haY{\widehat Y}

\def\bfPhi{\bfmath{\Phi}}

\newlength{\dhatheight}

\def\distarrow{\buildrel{\rm dist}\over\longrightarrow}

\def\thetaPR{\theta^{\text{\tiny\sf  PR}}}
\def\tilthetaPR{\tilde{\theta}^{\text{\tiny\sf  PR}}}

\def\zPR{Z^{\text{\tiny\sf  PR}}}  %Moved to upper case for consistency with M&T etc
 
\def\tTheta{{\text{\tiny$\Theta$}}}

\def\thetaPR{\theta^{\text{\tiny\sf  PR}}}
\def\tilthetaPR{\tilde{\theta}^{\text{\tiny\sf  PR}}}
\def\SigmaPR{\Sigma^{\text{\tiny\sf PR}}_\tTheta}

\usepackage{scalerel}
\newcommand{\smallcalS}{\ThisStyle{\scalebox{0.8}{$\SavedStyle\mathcal{S}$}}}

\def\EpsLength{\smallcalS}
\def\hafee{\widehat\fee}

\def\hac{\hat{c}}

\DeclareMathAccent{\widecheck}{0}{mathx}{"71}

\def\whamb{\wham{$\bullet$} }

\def\tTheta{{\text{\tiny$\Theta$}}}

\def\SigmaPR{\Sigma^{\text{\tiny\sf PR}}_\tTheta}

\def\SigmaTh{\Sigma_\tTheta}

\usepackage[usenames,dvipsnames]{color}
\usepackage{color}
\definecolor{programcode}{gray}{0.9}

\definecolor{lightgray}{gray}{0.7}

\definecolor{MyDarkBlue}{cmyk}{0.5,0.1,0,0.9}

\def\bl#1{{\color{blue}#1}}

\usepackage{cleveref}

\Crefname{corollary}{Corollary}{Corollaries}
\Crefname{eqnarray}{eq.}{eqs.}
\Crefname{equation}{eq.}{eqs.}

\Crefname{figure}{Fig.}{Figs.}
\Crefname{tabular}{Tab.}{Tabs.}
\Crefname{table}{Tab.}{Tabs.}
\Crefname{lemma}{Lemma}{Lemmas}

\Crefname{theorem}{Thm.}{Thms.}
\Crefname{definition}{Definition}{Definitions}
\Crefname{section}{Section}{Sections}
\Crefname{proposition}{Prop.}{Propositions}
\Crefname{assumption}{Assumption}{Assumptions}
\Crefname{example}{Example}{Examples}

\def\whamit#1{\smallbreak\pagebreak[3]%
\noindent\textit{#1}\ \ \gobblepars}

\def\wham#1{\smallbreak\pagebreak[3]%
\noindent\textbf{#1}\ \ \gobblepars}

\def\bl#1{{\color{blue}#1}}

\def\sstate{\textsf{S}}
\def\preDens{f_0}
\def\postDens{f_1}

\def\pFA{p_{\textsf{\tiny FA}}}

\def\Obs{Y}
\def\bfObs{\bm{Y}}

\def\postObs{X^1}
\def\bfpostObs{\bfmX^1}
\def\preObs{X^0}
\def\bfpreObs{\bfmX^0}

\def\condDist{\Uppi} %     Is this any better?
\def\InfoState{\clX}    
\def\surf{\breve{f}}
\def\surg{\breve{g}}

\def\thexp{\upupsilon}

\def\tchange{\uptau_{\text{\scriptsize\sf a}}}

\def\tstop{\uptau_{\text{\scriptsize\sf s}}}

\def\expa{\varrho_{\text{\scriptsize\sf a}}}

\newcommand{\overbar}[1]{\mkern 1.5mu\overline{\mkern-1.5mu#1\mkern-1.5mu}\mkern 1.5mu}

\def\thresh{%
\mathchoice
	{\text{\small\rm H}}%
	{\text{\small\rm H}}%
	{\text{\scriptsize\rm H}}%
	{\text{\tiny\rm H}}}

\def\barthresh{\overbar{\thresh}}

\def\MDD{{\sf MDD}}
\def\MDE{{\sf MDE}}

\def\ARL{{\sf ARL}}
\def\CauchyScale{\gamma}

\def\surL{\breve{L}}
\def\RegenState{\Delta}
\def\RegenTime{\tau}

\def\tregen{\tau_\RegenState}
\def\TD{D}
 
	\def\piX{\uppi}     %What was called \marg
\def\piPhi{\upvarpi}    %Now \varpi
 \def\muRegen{\upmu}    %For the regeneration theory 

\def\whamb{\wham{$\bullet$}}
\def\whamc{\wham{$\circ$}}
\def\barfinf{\barf_{\infty}}

\def\lips{\ell_{\text{\tiny\sf Z}}}         % For Lipschitz constants in Appendix to avoid using L. Sadly \ell used already in the form \ell_\phi.
\def\bigstate{\textsf{Z}}
\def\bz{\mathcal{B}(\bigstate)}

\def\modelAone{\text{1a}} % QCD Model cases
\def\modelAtwo{\text{1b}}
\def\modelAthree{\text{1c}}
\def\modelBone{\text{2a}}
\def\modelBtwo{\text{2b}}
\def\modelBthree{\text{2c}}
\def\modelCone{\text{3a}}
\def\modelCtwo{\text{3b}}

\def\bqed{{\color{blue} \qedb} \bigskip}

\title{Reinforcement Learning for Optimal Stopping in POMDPs with Application to Quickest Change Detection}

\author{Austin Cooper
\and Sean Meyn% <-this % stops a space 
\thanks{ASC and SPM are with the University of Florida, Gainesville, FL 32611 Financial support from ARO award W911NF2410389     and NSF awards  CCF 2306023 and DMS 2427265 is gratefully acknowledged.}%
}

\begin{document}
	
\maketitle

\begin{abstract}

The field of quickest change detection (QCD) focuses on the design and analysis of online algorithms that estimate the time at which a significant event occurs.   In this paper, design and analysis are cast in a Bayesian framework, where QCD is formulated as an optimal stopping problem with partial observations. 
An approximately optimal detection algorithm is sought using   techniques from reinforcement learning.   

The contributions of the paper are summarized as follows:

	\whamb
		A Q-learning algorithm is proposed for the general partially observed optimal stopping problem. It is shown to converge under linear function approximation, given suitable assumptions on the basis functions. An example is provided to demonstrate that these assumptions are necessary to ensure algorithmic stability.

	\whamb
		Prior theory motivates a particular choice of features in applying Q-learning to QCD. It is shown that, in several scenarios and under ideal conditions, the resulting class of policies contains one that is approximately optimal.
		
	\whamb
		Numerical experiments show that Q-learning consistently produces policies that perform close to the best achievable within the chosen function class.

\end{abstract}

\clearpage

\section{Introduction} 
\label{s:intro}

	The goal of the
	research surveyed in this paper is to create algorithms for \textit{quickest change detection} (QCD),
	for applications in which statistics are only partially known, particularly after the change has occurred.   
	Examples of events that we wish to detect include  human or robotic intruders, computer attack, faults in a power system,  and onset of heart attack~\cite{liatarvee21,liavee22}.

The standard QCD model includes a sequence of observations $\bfObs\eqdef \{ \Obs_k :  k\ge 0\}$ taking values in Euclidean space.   
	The statistics of these observations change at a time denoted $\tchange\ge 0$,  formalized through the representation 
\begin{equation}
	\Obs_k = \preObs_k \ind_{k< \tchange } +  \postObs_k \ind_{k \ge \tchange }  \,,\qquad k\ge 0\, .
	\label{e:QCDmodel}
\end{equation}    
The goal is to construct an estimate of the change time, denoted $\tstop$, that is adapted to the observations.  That is, on denoting $\Obs_0^k = (\Obs_0; \cdots; \Obs_k)$,  for each $k$ we may write $\ind\{\tstop \le k\} =  \fee_k(\Obs_0^k)$ for some Borel-measurable mapping $\fee_k\colon \ystate^{k+1} \to \{0,1\}$. We refer to any such sequence $\fee = \{ \fee_k : k \ge 0 \}$ as the \textit{policy}.
	The estimate $\tstop$ must balance  two costs: 
	1.~\textit{Delay}, which is expressed $(\tstop -\tchange)_+\eqdef \max(0,\tstop -\tchange)$,  and 2.~\textit{false alarm}, meaning that $\tstop -\tchange<0$.

	There are two general models that lead to practical solutions: minimax and Bayesian.
	The latter is the focus of this paper,   in which typical measures of performance are based on the classical Bayesian criteria of \textit{mean detection delay} $\MDD$ and probability of false alarm $\pFA$: 
\begin{equation}
	\text{ $\MDD  =   \Expect[ (\tstop -\tchange)_+ ]  $  
\ \ 
			and 
\ \ 
			$\pFA=  \Prob\{\tstop < \tchange \}  $.}
	\label{e:MDDpFA}
\end{equation}

	We opt for an optimality criterion that reflects the reality that there is little cost to trigger an alarm briefly before the change time.    
	The probability of false alarm is replaced by the
\textit{mean detection eagerness} $\MDE=\Expect[ (\tstop -\tchange)_- ]$.   For a fixed constant $\kappa>0$ and any policy $\fee$, denote    
\begin{equation}
		J(\fee)  = \MDD+\kappa\MDE =   \Expect\big[ (\tstop -\tchange)_+  +  \kappa (\tstop -\tchange)_-  \big].
	\label{e:MDD+kappaMDE}
\end{equation} 
Approaches to QCD with criterion \eqref{e:MDD+kappaMDE} are based on a partially observed Markov Decision Process (POMDP) model. 	A POMDP can be converted to a fully observed MDP through choice of state process, known as the  \textit{information state} \cite{put14}.    This leads to practical methods to obtain an optimal policy only for very simple models, such as Shiryaev's model introduced in \Cref{s:POMDPQCD}, in which $ \tchange$ has a geometric distribution.

Much of the literature on QCD is concerned with performance of heuristics,    typically based on the construction of a real-valued stochastic process $\{\InfoState_n\} $, and a policy defined by   a pre-assigned threshold $\thresh > 0$; the stopping rule is    
\begin{equation}
	\tstop = \min\{ n\ge 0 :   \InfoState_n\ge \thresh \}\, .
	\label{e:threshold}
\end{equation}%
\begin{subequations}%	
Two famous examples are  defined recursively:  with $\InfoState_0=0$,  
\begin{align}
			&\text{1. Shiryaev–Roberts:}   &&
\InfoState_{n+1} =   \exp\bigl( F_{n+1} \bigr) [\InfoState_n  +  1] 
\label{e:SR}
\\
			&\text{2. CUSUM:}   &&
\InfoState_{n+1} =    \max\{ 0,  \InfoState_n +  F_{n+1}  \} 
\label{e:CUSUM}
\end{align}
where $\{ F_{n+1} : n\ge 0\}$  is a stochastic process adapted to the observations. 

See the literature survey for highlights of the theoretical development of these tests.  

	\label{e:2surInfoState}
\end{subequations}

In the present paper,    statistics such as  	\eqref{e:2surInfoState} are proposed to define   a \textit{surrogate information state} (SIS) in the design of Q-learning algorithms.

\wham{Contributions}

\whamb
The paper begins with a treatment of Q-learning for the general optimal stopping problem.   
The need to incorporate partially observed cost prevents the application of standard algorithms such as those considered in \cite{tsivan99}.  
\Cref{t:Q} presents results for the total cost problem with discounting, while the more delicate case without discounting is addressed in \Cref{t:Qgamma1a}. The latter is most relevant to the QCD setting, since the objective \eqref{e:MDD+kappaMDE} does not involve discounting.

The proposed Q-learning algorithms are shown to converge under conditions similar to those imposed in \cite{tsivan99}, including under linear function approximation,  but further assumptions are required for the basis:   An example is provided to demonstrate that new assumptions on the basis are necessary to ensure algorithmic stability.
		
\whamb      
Theory for CUSUM surveyed in \Cref{s:asy} motivates the use of the statistic $\{\InfoState_{n+1}  \}$ defined in \eqref{e:CUSUM} as a SIS for Q-learning in applications to QCD.  While the observation models considered generalize those used in prior work on QCD (e.g., conditionally independent i.i.d.\ or Markov models), the features proposed for the Q-learning architecture do not satisfy the ergodicity conditions required in \Cref{t:Qgamma1a}. To address this, a regenerative implementation of Q-learning is introduced in \Cref{s:Qregeb} with convergence theory summarized in \Cref{t:Qgamma1}.

	\whamb
		Numerical experiments show that Q-learning consistently produces policies that perform close to the best achievable within the chosen function class. In some cases, the use of a \textit{multi-dimensional} SIS defined by several  CUSUM statistics   results in better performance than obtained with any of the threshold rules \eqref{e:threshold} using a single component of the  multidimensional statistic.

\wham{Literature}

	Our treatment of Q-learning in this paper for the general optimal stopping problem follows closely the framework of Tsitsiklis and Van Roy \cite{tsivan99}, where it is shown that Q-learning for optimal stopping converges under conditions similar to those for the simpler TD-learning algorithm.    This prior work considers only the discounted-cost setting,  and we are not aware of prior work addressing the challenges of partial observations.   
	
	Analysis of the threshold policy \eqref{e:threshold} is typically posed in an asymptotic setting, considering a sequence of models with threshold $\thresh$ tending to infinity.  Results establishing approximate optimality may be found in   \cite{mou86,mostar09} 
	 for either statistic \eqref{e:SR} or \eqref{e:CUSUM}.  
These results are restricted to the conditionally independent i.i.d.\ model in which each of the processes $\bfpreObs,\bfpostObs$ in \Cref{e:QCDmodel}  is i.i.d.,   and approximate optimality requires $F_{n+1} = L(Y_{n+1})$ in 
\eqref{e:CUSUM}, with $L$ the log likelihood ratio (LLR) between the marginal distributions of $\bfpreObs,\bfpostObs$.

While the early work of Shiryaev is cast in a Bayesian setting
\cite{shi77},  most of the literature is cast in a minimax setting in which $\tchange$ is modeled as unknown but deterministic and the metric of performance \eqref{e:MDDpFA} is replaced with Lorden's criteria: the \textit{worst case mean detection delay} $\sf{WADD} = \sup_n  \Expect[  (\tstop - n)_+ \mid \tchange = n]  $ and the \textit{average run length to false alarm},  $\ARL=\Expect[\tstop \mid \tchange=\infty]$.

Beyond  i.i.d.\ settings, analysis of CUSUM in \cite{lai98}  establishes lower bounds on $\sf{WADD}$ under the $\ARL$ constraint,  and analogous bounds on $\MDD$ under a constraint on $\pFA$ for the Bayesian regime.  Extensions to observation models that are conditionally Markovian or   hidden Markov models (HMMs) are found in \cite{zhasunherzou23,zhasunherzou23b}. 

 The potential \textit{sub-optimality} of CUSUM is a topic of   \cite{tarvee04} where it is shown that mean delay exceeds the optimal by a factor dependent on the geometric distribution for $\tchange$.

There is also a growing literature on
	robust QCD.  For example in 
\cite{unnveemey11,mouwudiadin24}
	it is assumed that the marginal distributions of 
	$\bfpreObs,\bfpostObs$ belong to distinct uncertainty classes, from which Least Favorable Distributions (LFDs) can be determined and incorporated into the CUSUM test. 
	For recent surveys on QCD theory see \cite{xiezouxievee21,liatarvee21,liavee22}.
	
The recent articles \cite{coomey24b,coomey24c}  obtain approximately optimal thresholds and approximations of the corresponding cost even in non-ideal settings for which $F_n$ is not associated with an LLR.

 	\smallskip

The present paper builds on the conference paper  \cite{coomey24d}, which is the first to propose the use of the CUSUM statistic in Q-learning for the QCD problem.  Convergence guarantees were absent, and we now understand that this lack of theory led to    Q-learning algorithms that were highly data inefficient.

\wham{Organization}  
	
\Cref{s:POMDP} describes a POMDP model for the general partially observed optimal stopping problem, along with the proposed  
Q-learning algorithm and   convergence theory.      Q-learning for the QCD problem is the subject of   \Cref{s:RLQCD}. 
 Numerical results illustrating the theory are contained  in \Cref{s:Qqcd}.	The Conclusions contain potential directions for future research.
 Some of the technical proofs and experimental details are contained in the Appendix.  
	
\section{Optimal stopping with partial observations}
\label{s:POMDP}
	
This section is entirely concerned with the general optimal stopping problem with partial observations.   
The results are applied to QCD in  \Cref{s:RLQCD}.
	
The partially observed Markov chain is denoted $\bfPhi$ and the observation process is denoted $\bfmY$,  evolving on respective sets $\state$,  $\ystate$.   These  are assumed to be Polish,  equipped with   Borel sigma-algebras $\bx$, $\by$ respectively (adopting the blanket assumptions of \cite{MT}).  It is assumed that $\Obs_k = h(\Phi_k)$,  $k\ge 0$,  for a measurable function $h\colon\state\to\ystate$.  
Note that i.i.d.\  measurement noise can be included by extending the definition of $\bfPhi$.  
In many of the results it is assumed that $\bfPhi$ has a unique invariant measure, denoted $\piPhi$.  

Many notational conventions from  the QCD model are maintained: 
The time we choose to stop is  denoted $\tstop$,   the input process takes on binary values, $U_k \in  \ustate = \{0,1\}$ and  $\tstop$ is defined as the first value of $k$ such that $U_k=1$. 
For given cost functions $c_\circ ,c_\bullet \colon\state    \to \Re$ and discount factor $0<\disc \le 1$, 
	the objective to be minimized over all   inputs adapted to the observations is
\begin{equation}
	\Expect\Big[ \sum_{k=0}^{\tstop -1} \disc^k c_\circ(\Phi_k  )     +   \disc^{\tstop}    c _\bullet(\Phi_{\tstop}  )    \Big]
	\label{e:ObjPOMDPQCD}
\end{equation}
This is a POMDP since $U_k$ is a function of present and past observations $Y_0^k$ for each $k$.  
	
Recall that a POMDP can be converted to a fully observed MDP through the introduction of an information state, so that an optimal policy is expressed as ``information state feedback''  \cite{put14,elllakmoo95}.  
The canonical choice is the sequence of  conditional distributions
of the state  $\Phi_k$  given the observations $Y_0^k$,  denoted $\{\condDist_k  :  k\ge 0\}$,
and the optimal policy is thus   $U^*_k = \fee^*(\condDist_k)$  for some binary-valued function $\fee^*$.

	To place the equation
	in standard form, denote $c(x,u) = (1-u) c_\circ(x) + u c_\bullet(x)$ for $x\in\state$ and $u\in \{0,1\}$.  
	In the POMDP model with  state process $\{\condDist_k  :  k\ge 0\}$,  the cost function
	is expressed as the conditional expectation $\clC(\condDist_k,U_k)    =   \Expect  \big[ c(\Phi_k, U_k)  \mid  Y_0^k \big] $, which is linear in $ \condDist_k$.    If 
	$\state$ is finite, with $K$ elements,  then $\condDist_k$ evolves on the $K$-dimensional simplex $\clS^K$,      
	so that  $\clC(\upbeta,u) = \sum_{x\in\state}  \upbeta(x) c(x,u)$ for $\upbeta\in \clS^K$.   
	
	With $Q^*_k  = Q^*(\condDist_k,U_k) $ and $\clC_k  = \clC(\condDist_k,U_k) $, the DP equation may be expressed
\[
	Q^*_k  = \clC_k
	+  \disc (1-U_k)  \Expect[ \uQ^*( \condDist_{k+1} )  \mid    Y_0^k\, ,  U_k  ]  
\]
where $\uQ(\condDist) = \min \{ Q(\condDist, 0), Q(\condDist, 1) \}$.

\smallskip

The formulation of Q-learning for optimal stopping begins with several design choices:

\wham{1.  Alternative to the information state}   Generation of the information state is complex, and requires a full model.    
 Rather, a \textit{surrogate information state} (SIS)   denoted $ \{\InfoState_n : n\ge 0 \}$  is constructed,  evolving on a subset of Euclidean space denoted $\sstate$, and adapted to the observations.   A common choice is a partial history of observations. 
 
 It is assumed without loss of generality that $\InfoState_n = g(\Phi_n)$,  $n\ge 0$, for some measurable function $g \colon\state\to\sstate$.
 
 \wham{2.   Function class}  
 Let   $\{Q^\theta (s,u) :  s\in\sstate\,, \ u\in\{0,1\}\,, \  \theta\in\Re^d\}$  denote a family of functions, continuously differentiable in $\theta$.        Associated with any $\theta\in\Re^d$ is the feedback law $\fee^\theta\colon\sstate\to\{0,1\}$ which defines the stopping rule,
\begin{equation}
	\fee^\theta(s) =  \ind\{ Q^\theta(s, 0) \ge Q^\theta(s, 1) \}     \,,  \qquad s\in \sstate \, .
	\label{e:fee_theta_HTorQCD}
\end{equation}

 \wham{3.   Criterion of fit}    The parameter $\theta^*\in\Re^d$ is optimal if $Q^{\theta^*}$ is a solution to   the \textit{projected Bellman equation}  \cite{sutbar18,CSRL}.  This may be expressed as the root finding problem  $\barf(\theta^*) =0$, with
\begin{equation}
\begin{aligned}
\barf(\theta) & \eqdef 
\Expect\big[  \zeta^\theta_n \{-  Q^\theta_n  +  \clC_n
			+  \disc (1-U_n)    \uQ^\theta(\InfoState_{n+1} )   \}  \big]
\\
			&=  \Expect [  \zeta^\theta_n  \TD^\theta_{n+1} ]\,,  
\end{aligned}
	\label{e:pBe}
\end{equation}
where  $  \zeta^\theta_n  $ is a function of $(\InfoState_n, U_n)$  [a typical choice is 
  $\zeta^\theta(s,u) = \nabla_\theta Q^\theta (s,u)$], 
	and
	$  \TD^\theta_{n+1} \eqdef  -  Q^\theta_n  +   c(\Phi_n, U_n)  
	+ \disc (1-U_n)    \uQ^\theta(\InfoState_{n+1} )  $ is the (fixed parameter) temporal difference.
	The second equation, 
	in which $ \clC_n$ is replaced by $c(\Phi_n, U_n)  $, 
	follows from the smoothing property of conditional expectation.  
	
The expectation in \eqref{e:pBe} is taken in ``steady-state'', whose definition requires assumptions on both the hidden state process $\bfPhi$ and the input sequence, which brings us to the fourth design choice. 

\wham{4.  Exploration}   
Randomized policies are traditionally required to ensure convergence, and are often parameter dependent. 
  We adopt the notation $\feex^{\theta}$  which defines the input as follows:  
\begin{equation}
\Prob\{ U_n = u\mid Y_i, U_i : i<n;  \   \InfoState_n = x   \}    
=
\feex^{\theta_n} (u \mid x)    
\label{e:feex}
\end{equation} 
It is shown in \cite{mey24} that the standard Q-learning  algorithm with linear function approximation is stable if $\feex^{\theta}$  is defined to approximate the $Q^\theta$-greedy policy $	\fee^\theta$.

 \wham{Q-learning and stochastic approximation}
 The root finding problem  $\barf(\theta^*) =0$ with $\barf$ defined in   \eqref{e:pBe}  motivates the   update rule,
\begin{equation}
	\theta_{n+1} =  \theta_n + \alpha_{n+1}  \zeta_n  \TD_{n+1}
	\label{e:POMDPQ}
\end{equation}
in which  $  \zeta_n  =  \zeta_n^{\theta_n}$,  $  \TD_{n+1}  = \TD_{n+1}^{\theta_n}$, and
the step-size sequence $\{ \alpha_{n+1}  \}$  is non-negative.    The input sequence $\{U_n\}$  used in the algorithm is assumed to be defined by a  
randomized stationary policy as in \eqref{e:feex}.

The recursion \eqref{e:POMDPQ} will be cast as an instance of stochastic approximation (SA).  Our main results are cast in SA theory, which requires the introduction of terminology from these field.

Recent results in \cite{borchedevkonmey25} consider the general SA recursion  $\theta_{n+1} =  \theta_n + \alpha_{n+1}  f(\theta_n, \qsaprobe_{n+1})$ in which $\bfqsaprobe$ evolves on some Polish space $\bigstate$, and $f \colon \Re^d \times \bigstate\to\Re^d$ is Borel measurable.  Additional assumptions on $f$ and $\bfqsaprobe$ imply convergence of $\{\theta_n \}$ (both almost surely and in mean square),   along with moment bounds and a Central Limit Theorem.

The starting point of analysis is the definition of the \textit{mean flow vector field}  $\barf(\theta) = \Expect[ f(\theta,\qsaprobe_k)]$,  $\theta\in\Re^d$,     where the expectation is in steady-state in a stationary realization of $\bfqsaprobe$;   this definition is consistent with 	\eqref{e:pBe}.
The SA recursion may be motivated by the associated   \textit{mean flow},
\begin{equation}
\ddt \odestate = \barf(\odestate) \,, \quad \odestate_0\in\Re^d\,,
\label{e:meanflow}
\end{equation}
which is typically assumed to be globally asymptotically stable to some stationary point $\theta^*\in\Re^d$.  
The strongest conclusions for the SA recursion  are obtained if  the mean flow is \textit{exponentially asymptotically stable} (EAS):  
 There is $b_0>0$, $\epsilon_0>0$ such that  for any initial condition and $t\ge0$,
    \begin{equation}
    	\| \odestate_t - \theta^* \|  \le  b_0 \exp(-\epsilon_0  t)   \| \odestate_0 - \theta^* \|.
    	\label{e:EAS}
    \end{equation}

A key assumption for the Markov chain is a Lyapunov drift condition  known as (DV3):
for functions $V\colon\bigstate\to\Re_+$,  $ W\colon\bigstate\to [1, \infty)$  and $b>0$,   
\begin{equation}
	\Expect\bigl[  \exp\bigr(  V(\qsaprobe_{k+1})      \bigr) \mid \qsaprobe_k=x \bigr]  
	\le  \exp\bigr(  V(x)  - W(x) +  b   \bigl) \,, \qquad x\in\bigstate.
	\label{e:DV3body}
\end{equation}
This is slightly weaker than the condition appearing in \cite{borchedevkonmey25},  but is equivalent subject to the assumptions on $W$ imposed below (in particular, that the sublevel sets of $W$ are \textit{small}---we defer to  the reference and \cite{MT} for  Markov chain terminology).    Markov chains that are uniformly ergodic satisfy the assumptions of  \cite{borchedevkonmey25},  which includes ergodic Markov chains on a finite state space.

We say that  a function $G\colon\bigstate \to\Re $ satisfies the growth condition $G = o(W)$ if 
\[
\lim_{r\to\infty}  \sup_{x\in\bigstate}  \frac{ |G(x)|}{\max\{r,W(x)\}} =0.  
\]
In some results of  \cite{borchedevkonmey25} it is assumed that $G_\theta = o(W)$ for each $\theta$,   where $G_\theta(\varble) =  \| f(\theta,\varble) \| $.

    %%%%%%%%%%%%%%%  Returning to Q-learning 

\smallskip

Given the form of the Q-learning recursion \eqref{e:POMDPQ}, to frame analysis within the setting of  \cite{borchedevkonmey25} it is assumed henceforth that $\qsaprobe_{k+1} = (\Phi_{k+1}; \Phi_k; U_k  )$ defines a Markov chain with state space   $\bigstate = \state\times\state\times \ustate $.

We also restrict to a linear function class, so that $Q^\theta = \theta^\transpose \psi$, and adopt the notation
\begin{equation}
	\psi(s,u) = u \psi^1(s) + (1-u) \psi^0(s)     \,, \qquad    s\in\sstate\,,  \ u\in \ustate
	\label{e:psiQ}
\end{equation}
where   $\psi^i : \sstate \to \Re^d $ are Borel measurable,  and to save space denote  $\psi_{(k)} \eqdef \psi(\InfoState_k, U_k)$  for each $k$.   
The vector field   \eqref{e:pBe}   may be expressed,   
\begin{equation}
	\barf(\theta)  =  
		- R \theta + b + p _0  \disc   \Expect\big[ \psi^0_{(k)} \uQ^\theta(\InfoState_{k+1} ) \big] 
\label{e:barfQ}
\end{equation}
where  $R \eqdef \Expect_\piPhi[ \psi_{(k)} \{\psi_{(k)}\}^\transpose]$,
$p_i=\Prob\{ U_k = i \}$ for $i=0,1$ (assumed independent of $k$), 
and
  $\uQ^\theta(s)=\min\{ \theta^\transpose \psi^{0}(s),  \theta^\transpose \psi^{1}(s)\}$.

We say the basis is \textit{separable}
 if 
$\psi^0_i (s)   \psi^1_i(s) = 0$ for every $s\in\sstate$ and $i$.
This assumption is imposed to establish EAS of the mean flow with vector field \eqref{e:barfQ}.
It is assumed 
  without loss of generality
  that the entries are ordered so that for some $0< d_0 <d$ we have $\psi^1_i \equiv 0$ for $i\le d_0$ and $\psi^0_i \equiv 0$ for $i> d_0$.
This justifies the notation
$\psi(s,u) = [ (1-u) \psi^{(0)}(s)  ;  u  \psi^{(1)}(s) ]$ in which  
	$\psi^{(j)}  \colon\sstate \to \Re^{d_j},  \, j=0,1$,
	are defined by $\psi^{(0)}  =  \psi^0_i $ for $1\le i \le d_0$, and $\psi^{(1)}_i  = \psi^1_{i-d_0} $ for $1\le i\le d_1 = d-d_0$.

	Write $\theta = [ \theta^0;\theta^1 ]$ with $\theta^0 \in \Re^{d_0}  $,  $\theta^1 \in \Re^{d_1}  $. 
	The proposition that follows provides conditions for convergence of the Q-learning recursion \eqref{e:POMDPQ} to some $\theta^*=(\theta^{*0};\theta^{*1})\in\Re^d$.  
	The representation \eqref{e:ProjectedCost} implies that 
	$ \theta^{*1}$ defines the projection of the cost $c(\varble, 1)$  onto the span of $\{\psi^{(1)}_i  : 1 \le i \le d_1\}$ within the Hilbert space $L_2(\piPhi)$.   
	This result imposes $\disc<1$, as in  \cite{tsivan99}.  Extension to $\disc =1$ is postponed to \Cref{t:Qgamma1a}.

\begin{proposition}[Conditions for convergence]
	\label[proposition]{t:Q}  
		Consider the Q-learning recursion \eqref{e:POMDPQ} with linear function approximation, subject to the following additional assumptions:
		
	\whamc  The step-size sequence is 
		$\alpha_n =  \min\{ \alpha_0  , n^{-\rho} \}$ with $\alpha_0>0$ and $1/2 < \rho \le 1$.

	\whamc
		The input $\{U_k : k\ge 0 \}$ used for training is i.i.d.\ on $\ustate = \{0,1\}$ and independent of $\bfPhi$. 
		
	\whamc The drift condition  \eqref{e:DV3body}
		holds for the Markov chain $\bfPhi$.  		
		Moreover, 
		$ S_W(r)  \eqdef \{ x \in\state :  W(x)\le r \}  $     is either small or empty and $ 
	\sup\{ V(x) :  x\in S_W(r) \}  <\infty$ for each $r>0$.

\whamc   $G = o(W)$, where  
$G(x) =   \max_u  [  \| \psi (g(x),u)  \|^2   + c(x,u)^2  ]$  \emph{(recall that $\InfoState_n = g(\Phi_n)$).}

\whamc  
		The basis is separable,  and $R \eqdef \Expect_\piPhi[ \psi_{(k)} \{\psi_{(k)}\}^\transpose]$  has full rank.
		
	\smallskip

	\noindent
		Then, 
	\whamrm{(a)}  For each initial condition $\theta_0,\Phi_0$,  the Q-learning recursion 
	\eqref{e:POMDPQ}  is convergent   to  some $\theta^*\in\Re^d$:
	\[
	\lim_{n\to\infty} \theta_n  = \theta^* \ \  a.s.\
	\]

	\whamrm{(b)}   The second component of $\theta^*$ may be expressed 
\begin{equation}
\theta^{*1}   =  R_{(1)}^{-1}  \Expect_\piPhi[  \psi_{(k)}^{(1)}  c _\bullet (\Phi_k) ]     \,,
\ \  \textit{	with $R_{(1)} = \Expect_\piPhi[  \psi_{(k)}^{(1)}  \{\psi_{(k)}^{(1)}\}^\transpose]$.    } 
\label{e:ProjectedCost}
\end{equation}
		
The function $\hac _\bullet  \eqdef  ( \theta^{* 1 } )^\transpose \psi^{(1)}$ defines the projection of $c_\bullet(\Phi_k)$ onto the span of $\{\psi^{(1)}_i \}$ in $L_2(\piPhi)$.  

	\whamrm{(c)}  Consider the estimates obtained via  averaging:   
\begin{equation}
\thetaPR_n = \frac{1}{n}\sum_{k=1}^n \theta_k \,,\qquad n\ge 1.
\label{e:thetaPR}
\end{equation}
Suppose that the step-size
$ \alpha_{n + 1}  =1/ (n+1)^\rho$ is used  with $\half < \rho < 1$.  Then, the CLT holds for the normalized sequence 
$\zPR_n\eqdef  \sqrt{n} \tilthetaPR_n$,
with $\tilthetaPR_n = \thetaPR_n -\theta^*$.   
Moreover,  the rate of convergence is optimal,
in that,
\begin{equation}
	\lim_{n \to \infty} n \Expect [ \tilthetaPR_n (\tilthetaPR_n)^\transpose  ] 
		= 	\lim_{n \to \infty}  \Expect [ \zPR_n (\zPR_n)^\transpose  ] = \SigmaPR \,,
\label{e:PR}
\end{equation}
where the covariance matrix $\SigmaPR$  is minimal in a matricial sense, achieving the lower bound of Chung-Polyak-Ruppert \cite{chedevbusmey20b}. 
\qed
\end{proposition}
\smallskip

Proof of the proposition may be found in \Cref{s:SAthy}. A key step is given below:

\begin{lemma}
	\label[lemma]{t:QstableMeanFlow}
	Under the assumptions of \Cref{t:Q}  the mean flow is  EAS.   
\end{lemma}

\wham{Proof}  
We first establish the existence of a unique solution to $\barf(\theta^*) = 0$.
Steps in this proof lead to the following inequality:  for some $\epsy_0>0$ and each initial condition for the mean flow, $\ddt \| \odestate_t - \theta^* \|^2  \le - \epsy_0   \| \odestate_t - \theta^* \|^2$.    This implies EAS.

	Separability combined with the independence assumption on $\{ U_k \}$ implies that the autocorrelation matrix admits a two-block diagonal form:
\[
R    =    \begin{bmatrix}  p_0 R_{(0)}  &  \Zero  \\  \Zero  & p_1  R_{(1)} \end{bmatrix}  
\]
in which $p_i = \Prob_\piPhi\{ U_k = i \}$,  and  $R_{(i)} = \Expect_\piPhi\big[ \psi_{(k)}^{(i)} {\psi_{(k)}^{(i)}}^\transpose\big]$ for each $i=0,1$.
The full rank assumption on $R$ implies that $p_i>0$ for each $i$, and that each  $R_{(i)} $ is positive definite. 

Separability implies that the mean flow dynamics for $\theta^1$ are linear: 
\begin{align*}
	\ddt \odestate^1_t &= p_1 \Expect_\piPhi[\psi_{(k)}^{(1)} \{ -\theta^{1 \transpose}\psi_{(k)}^{(1)} + c _\bullet (\Phi_k) \} ]
	\\
	&= - p_1 R_{(1)} \theta^1 + p_1 \Expect_\piPhi[  \psi_{(k)}^{(1)}  c _\bullet (\Phi_k) ]
\end{align*}
The matrix $-p_1 R_{(1)} $ is Hurwitz, and hence $\theta^1_t$  is convergent to the value  $\theta^{*1}$ given in \eqref{e:ProjectedCost}.

To establish convergence of $ \odestate^0_t$,  we introduce $H^{\theta^0}(s) = Q^\theta(s,0) =\theta^{0 \transpose}\psi^{(0)}(s) 
$, which is the function of interest in   \cite{tsivan99}.
From \Cref{e:pBe} we have,
\begin{equation}
	\ddt \odestate^0_t = \barf^0(\odestate_t) \,, \qquad 
	\barf^0(\theta)=  p_0 \Expect[ \psi^{0}_k D^\theta_{k+1}]
	\label{e:tsivan99-1}
\end{equation} 
where $D^\theta_{k+1} = -  H^{\theta^0}(\InfoState_k)  + c_\circ(\Phi_k) + \disc \uQ^\theta(\InfoState_{k+1})$ with $\uQ^\theta(s) = \min \{ H^{\theta^0}(s)   , \theta^{1 \transpose}\psi^{(1)} (s) \}$.    
Since we have established convergence of $\theta^1_t$,    solutions to \eqref{e:tsivan99-1}
converge to those of the reduced order ODE $\ddt x_t = \barf^0 (x_t  ,  \odestate^{*1}) $, which is equivalently expressed
\[
\begin{aligned}
	&\ddt x_t = 
	p_0 \barg( x_t )\,,
	\quad 
	\barg(x) 
	= \Expect[  \psi_{(k)}^{(0)}  D^x_{k+1}   ]  \,,
	\\
	&  D^x_{k+1} \hspace{-.2em} = \hspace{-.2em} - H^x   (\InfoState_k)   +  c_\circ(\Phi_k) + \disc \min\{ H^{ x}   (\InfoState_{k+1}) ,  \hac _\bullet (\InfoState_{k+1}) \} 
\end{aligned}
\] 
where  $\hac _\bullet$ is defined in \Cref{t:Q}. 

The vector field $\barg$ coincides with the mean flow obtained in  the fully observed optimal stopping problem of   \cite[Section III-C]{tsivan99}.  
Application of the main result there  establishes convergence of $x_t$;     we require  $R_{(0)} >0 $ and $p_0>0$, which holds 
under our assumption that $R>0$.
\qed

\begin{proposition}[Potential for algorithmic instability]
\label[proposition]{t:Qbad}   
There is an example satisfying  all assumptions of \Cref{t:Q} except separability,   such that the  mean flow is not stable:  There is a unique root $\theta^*$ to $\barf(\theta^*) =0$,  yet for any $\theta \neq \theta^*$ the solution to the mean flow   \eqref{e:meanflow}  diverges to infinity for $\odestate_0\neq \theta^*$.
\end{proposition}

\wham{Proof}   The counterexample is constructed through the following steps:
\wham{1.}   $\psi^1 = \xi \psi^0$.    Moreover,   $R_{(0)} 
= \Expect_\piPhi[\psi^0_{(k)} (\psi^0_{(k)})^\transpose]$ is full rank,
and there is a positive definite matrix $R_{(0)}^-$ such that
$ \Expect_\piPhi[ ( \{ \theta^\transpose \psi^0_{(k)} \}_- )^2  ] \ge \theta^\transpose  R_{(0)}^-  \theta
$  for all $\theta\in\Re^d$, 
where $z_- = \max(0,-z)$ and  $z_+ = \max(0,z)$ for $z\in\Re$.

As just one example, 
suppose that $\{  \psi_{(k)}^0    : k\ge 0\}  \eqdist \{ - \psi_{(k)}^0    : k\ge 0\}  $  (the finite dimensional distributions coincide for a steady-state realization).   In this case we may take $ R_{(0)}^- = \half R_{(0)}$.

\wham{2.}   $\delta = \Prob_\piPhi[ \psi^0_{(k+1)}  \neq \psi^0_{(k)} ] $ is vanishing,  with $\delta(\xi) \xi = o(1)$   as $\xi\to\infty$.

\wham{3.}  
For the i.i.d.\ input with $p_0 = \Prob\{U_k=0\}$,  $p_1=1-p_0$,    we have $ p_1(\xi) \xi^2 = o(1)$ as $\xi\to\infty$.

\wham{4.}   The vector $b$ appearing in  \eqref{e:barfQ} is zero,   $b =  \Expect_\piPhi[\psi_{(k)} c(\Phi_k,U_k)]$.   Hence $\theta^* =0$ solves $\barf(\theta^*) =0$.

The vector field \eqref{e:barfQ} admits the approximation,
\[
\begin{aligned}
	\barf(\theta)  &= -(p_0 + p_1 \xi ) R_{(0)}  \theta   
	+  p_0 \disc \Expect \big[  \psi_{(k)}^0  \min \{  \theta^\transpose \psi_{(k+1)} ,  \xi \theta^\transpose \psi_{(k+1)}  \}  \big]  
	\\
	&=
	- R_{(0)}  \theta   
	+    \disc \Expect \big[\psi_{(k)}^0  \min \{  \theta^\transpose \psi^0_{(k)} ,  \xi \theta^\transpose \psi^0_{(k)}  \}  \big]  
	+  \clE (\theta)  
\end{aligned}
\]
in which  $\| \clE (\theta)  \| \le   \epsy_\xi  \|\theta\|$ for all $\theta,\xi$,  and $\epsy_\xi =o(1)$.   
Note that (2) was applied here to replace $ \theta^\transpose \psi_{(k+1)}$  by  $ \theta^\transpose \psi_{(k)}$  in the minimum,
and (1) was then used to replace $   \psi_{(k)}$ by $   \psi_{(k)}^0$.

Under the assumption that $\xi>1$ we obtain
\[
\begin{aligned}
	\theta^\transpose \barf(\theta)  
	& =  -  \theta^\transpose R_{(0)}  \theta    
	+    \disc \Expect \big[ \theta^\transpose \psi_{(k)}^0   \min \{  \theta^\transpose \psi^0_{(k)} ,  \xi \theta^\transpose \psi^0_{(k)}  \}  \big]  
	+ \theta^\transpose \clE (\theta)  
	\\
	& = 
	-  \theta^\transpose R_{(0)}  \theta +
	\disc   \Expect \big[  (  \{ \theta^\transpose \psi_{(k)}^0  \}_+ )^2  + \xi(  \{ \theta^\transpose \psi_{(k)}^0  \}_- )^2     \big]   + \theta^\transpose \clE (\theta)  
	\\
	& \ge    - \theta^\transpose R_{(0)}  \theta    +   \disc  \xi   \theta^\transpose R_{(0)}^-  \theta   -  
	\| \clE (\theta)  \|   \|\theta\| 
\end{aligned}
\]

From this we conclude that there is $\nu>0$ such that for all $\xi>1$ sufficiently large and all $\theta$ we have the lower  bound $\theta^\transpose \barf(\theta)   \ge \nu \|\theta\|^2 $.
It follows that $\ddt \| \odestate_t \|^2  \ge 2  \nu \|\odestate_t\|^2$ for each initial condition, which implies the desired conclusion.   
\qed

\medskip
 
We now consider conditions for convergence without discounting.  
In the following extension of \Cref{t:Q}  
  we require that the  \textit{conditional} covariance is full rank, denoted
\begin{equation}
	\Sigma_{k+1\mid k} =
	\Expect_\piPhi[  ( \psi_{(k+1)} - \psi_{(k+1 \mid k)}  ) ( \psi_{(k+1)} - \psi_{(k+1 \mid k)}  ) ^\transpose]
	%\Expect_\piPhi[  (\psi_{(k+1 \mid k)}  - \psi_{(k+1)})   (\psi_{(k+1 \mid k)}  - \psi_{(k+1)})  ^\transpose] 
	\label{e:condCov}
\end{equation}
	with 
	$\psi_{(k+1 \mid k)} =  \Expect[\psi(\InfoState_{k+1}, U_{k+1})\mid \Phi_k]$.

\begin{proposition}[Convergence without discounting]
	\label[proposition]{t:Qgamma1a}   
		Suppose the assumptions of \Cref{t:Q} hold and in addition $\Sigma_{k+1\mid k} $ has full rank.   Then the conclusions of 
	\Cref{t:Q} hold for $\disc =1$.       
\end{proposition}
	
\wham{Proof}  
	The proof of the proposition requires a closer examination of the convergence proof in \cite{tsivan99}.    Under the assumptions of the proposition we establish that the convergence proof of  \cite{tsivan99} extends to $\disc =1$,  and the remainder of the proof follows that of  \Cref{t:Q}.

Subject to the assumptions of the proposition,   we show that $\disc    P$  is a strict contraction in the 
	subspace of $L_2(\piPhi)$ spanned by $\psi$,  denoted $\clH =\{  \theta^\transpose \psi : \theta\in\Re^d \}$.     
	The contraction property holds even when $\disc=1$:   there is $\varrho <1$ such that $\| P g \|_\piPhi^2  \le \varrho   \| g\|_\piPhi^2  $  for $g\in \clH$.     This contraction property is all that is required in  \cite{tsivan99}.

Establishing the contraction begins with the definitions 
\begin{equation}
\begin{aligned} 
\| P g \|_\piPhi^2   &=   \Expect_\piPhi[   \Expect[g(\Phi_{k+1})\mid \Phi_k] ^2 ]  
\\
\|  g \|_\piPhi^2   &=  \Expect_\piPhi[   g(\Phi_{k+1}) ^2 ]  =
 \| P g \|_\piPhi^2  +  \sigma^g_{k+1\mid k} 
\end{aligned} 
	\label{e:Contract_a}
\end{equation}
where $\sigma^g_{k+1\mid k} = \Expect_\piPhi[  \{ g(\Phi_{k+1}) -  \Expect[ g(\Phi_{k+1})\mid \Phi_k] \}^2 ]  $.

For $g= \theta^\transpose \psi$ we have  $ \| g \|_\piPhi^2 = \theta^\transpose R \theta$ and   applying 	\eqref{e:Contract_a},
\[
\| P g \|_\piPhi^2  =  \|  g \|_\piPhi^2   -  \sigma^g_{k+1\mid k}     =  \|  g \|_\piPhi^2   -     \theta^\transpose   \Sigma_{k+1\mid k}        \theta  
\]	
By assumption we   have 	$  \Sigma_{k+1\mid k}  \ge (1-\varrho) R$  for some $\varrho>0$, and hence
$\| P g \|_\piPhi^2   \le  \varrho  \| g \|_\piPhi^2 $.
\qed

\section{Bayesian QCD and RL}
\label{s:RLQCD}

	The theory in the previous section is applied here to the   Bayesian QCD problem.    For this we require the specification of a suitable POMDP,  along with the choice of SIS (surrogate information state) in a Q-learning architecture---recall the discussion at the start of \Cref{s:POMDP}.   The brief survey in \Cref{s:asy} is intended to motivate the use of the CUSUM statistic as a SIS.

\subsection{POMDP model for QCD}
\label{s:POMDPQCD}
	
	The notation remains the same as in \Cref{s:POMDP}:  the hidden Markovian state process   $\bfPhi$ and observation process   $\bfmY$ evolve on Polish spaces $\state$,  $\ystate$.   The observations  are a deterministic function of the state process: $Y_k = h(\Phi_k)$   for a measurable function $h\colon\state\to\ystate$.  
	
We henceforth use \textit{QCD-POMDP} to refer to the model whose special structure is   described as follows:
	
\whamc     There is a decomposition $\state = \state_0 \cup \state_1$,   for which $\state_1\in\bx$ is \textit{absorbing}:   $\Phi_k\in \state_1$ for all $k\ge 0$ if $\Phi_0\in \state_1$.      The change time is defined by $\tchange = \min\{k\ge 0 :  \Phi_k \in \state_1 \}$.  
	
\whamc    In the objective function \eqref{e:ObjPOMDPQCD} we take  
\[
		c_\circ(x)   =  \ind\{ x \in \state_1\}    \,, \quad
		c_\bullet(x) = \kappa  \Expect[ \tchange \mid \Phi_0=x] \ind\{ x\in\state_0 \}
\]
	which is consistent with  \eqref{e:MDD+kappaMDE}   when  $\disc =1$.  	
	
This POMDP model captures many change detection models considered in past literature.

\wham{Preventative maintenance}
	The choice of cost criterion \eqref{e:MDD+kappaMDE} with mean detection eagerness $\MDE=\Expect[ (\tstop -\tchange)_- ]$  is chosen in part for application to  optimal decision making for preventative maintenance, which concerns problems similar to QCD, often in a Bayesian setting \cite{dejsca20,mcc65}.   The cost criteria are slightly different, including  metrics such as the mean time between system failures and mean time to repair. Much of the literature adopts a POMDP model that is a special case of the one proposed here  \cite{morpapandnierig22}.

\wham{Independent model}
Also called the \textit{conditionally} independent model,
this is  defined by  \eqref{e:QCDmodel} in which  $\bfpreObs$ and $\bfpostObs$ are two  mutually independent stochastic processes on the common state $\ystate$ that are also independent of the change time $\tchange$.   The most common examples impose additional assumptions: 
\wham{$\circ$}  
Independent i.i.d.\     model, in which  $\bfpreObs$ and $\bfpostObs$  are each i.i.d.\ sequences.

A special case is 
\textit{Shiryaev's model}, 
 in which   $\tchange$ has a geometric distribution.   This provides an example  of the QCD-POMDP model in which $\bfPhi =  ( \bfpreObs;  \bfpostObs;  \bfmI)$  with $I_k = \ind\{ \tchange  \le k \}$,  and  $\state = \ystate\times \ystate \times \{ 0, 1\}$.

 \wham{$\circ$}  
Independent Markovian model, in which   $\bfpreObs$ and $\bfpostObs$  are each time-homogeneous Markov chains.

If in addition the change  time $\tchange$ is defined as the sojourn time to some set for a Markov chain $\bfmW$, 
then we obtain an instance of the  QCD-POMDP model with $\bfPhi = (\bfpreObs ;\bfpostObs ;\bfmW)$.   
If $\bfPhi$ evolves   on a finite state space then  an optimal policy may be obtained numerically through the introduction of an  information state.   	In Shiryaev's model, we may take    the real-valued process $\{ \beta_k =  \Prob\{  \tchange \le  k \mid \Obs_0, \dots, \Obs_k \}  : k\ge 0\}$, and an optimal test is of the form $U_k^* = \ind\{  \beta_k \ge \thresh\} $ for some   threshold $\thresh>0$  \cite{shi77,veeban14}.

	Outside of such special cases, computation and implementation of an optimal test is  far too complex to be practical.

\subsection{Asymptotic theory}
\label{s:asy}

We survey here theory from \cite{coomey24b,coomey24c} concerning approximately optimal tests based on the adoption of the CUSUM statistic \eqref{e:CUSUM} as a SIS.   This will motivate Q-learning architectures described in \Cref{s:Qqcd}.     While this prior work develops theory for the general 
QCD-POMDP model, the survey here is limited to the independent i.i.d.\ model for which 
 stochastic processes $\bfpreObs,\bfpostObs$ in \eqref{e:QCDmodel}
 are i.i.d., with marginal distributions denoted $\piX^0,\piX^1$ respectively.

The approximations are based on large $\kappa>0$ for the criterion \eqref{e:MDD+kappaMDE}. 
Analysis hinges on mild assumptions on the  function $F$ for CUSUM  (the function $F$ defines $F_n=F(Y_n)$ in \eqref{e:CUSUM}),  and also a regularity condition for the change time.     Denote the log moment generating functions $\Lambda_i(\thexp)=\log\Expect[\exp(\thexp F(X^i_k))]$,   $\thexp\in\Re$,  $i=0,1$.

\wham{(A1)} \textit{Moment conditions}: 
The function $F\colon\ystate\to\Re$ is Borel measurable. Letting $m_i=\int F(s)\,f_i(s)\, ds$ for $i=0,1$, it is assumed that $m_0<0$ and $m_1>0$.

\wham{(A2)} \textit{Regular geometric tail}: for some $\expa\in(0,\infty)$,
\begin{equation}
\lim_{n\to\infty}\frac{1}{n}\log\Prob\{\tchange\ge n\}=-\expa 
\label{e:hazardAss}
\end{equation}

\wham{(A3)}  
  The two log moment generating functions are finite valued in a neighborhood of the origin.   In addition,  $\Lambda_0$ has two distinct roots $\{0,\thexp_0\}$, a unique solution $\thexp_+>\thexp_0$ to $\Lambda_0(\thexp_+)=\expa$, and is finite in a neighborhood of $[0,\thexp_+]$.  

\smallskip
Under (A1), it follows that the CUSUM statistic $\{\InfoState_n\}$ evolves as a reflected random walk (RRW) with negative drift for $n<\tchange$ and thereafter a RRW with positive drift.    The regularity assumption (A2) obviously holds in Shiryaev's model for which $ n^{-1}\log\Prob\{\tchange\ge n\}$ is independent of $n$.   

 Suppose that the marginals  $\piX^0,\piX^1$  are mutually absolutely continuous with LLR denoted $L=\log(d\piX^1/d\piX^0)$. 
This function satisfies the sign conditions of (A1):
\begin{equation}
m_0 = \piX^0(L)=-D(\piX^0\|\piX^1)<0 \,, \ \   m_1 = \piX^1(L)=D(\piX^1\|\piX^0)>0
\label{e:signcon1}
\end{equation}
where $D$ denotes relative entropy. It is known that the use of $F_n=L(Y_n)$ in \eqref{e:CUSUM} defines a test that is approximately optimal under the traditional Bayesian cost criteria \eqref{e:MDDpFA}.  An extension of this theory is summarized in \Cref{t:barcApprox}
below.

\begin{figure}[t]
	\centering		
	\includegraphics[width=0.6\hsize]{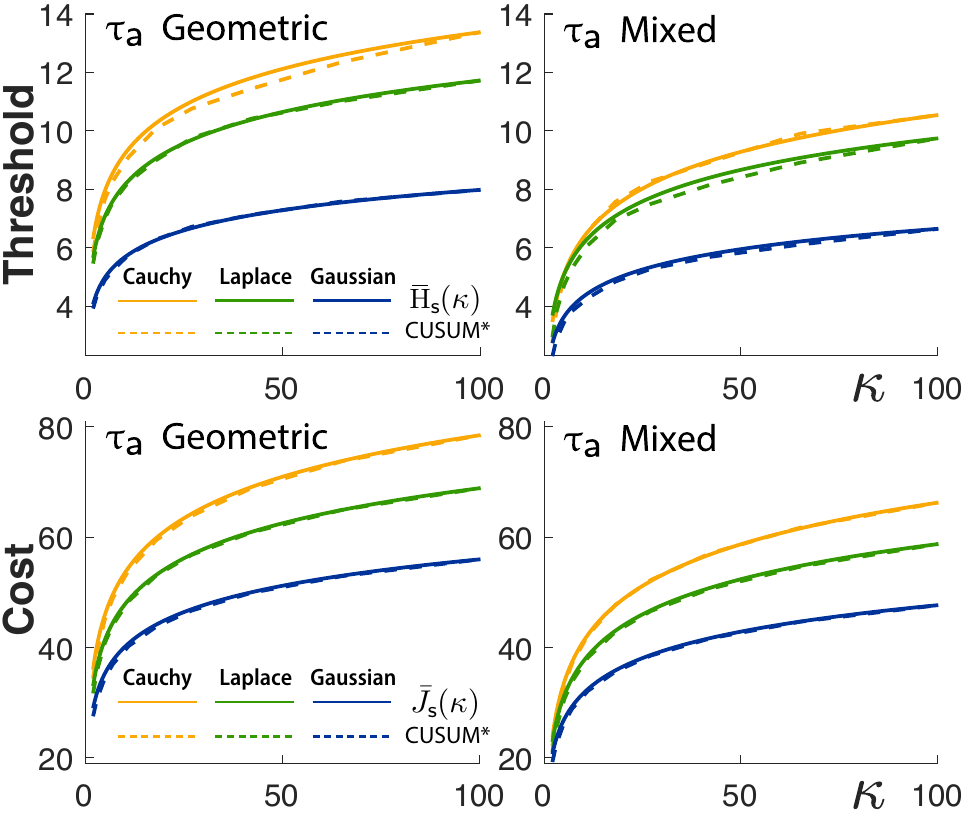} %trying to keep this fig on current page
	\caption{Approximations \eqref{s:shiftHJ} of the optimal threshold and cost compared to CUSUM*. }
	\label[figure]{f:CostAndThresholdApproximations}
\end{figure}
	
\begin{subequations}

Once we have chosen the stochastic process $\{ 	F_{n+1}  \} $ in the CUSUM statistic 
$\{ \InfoState_{n+1} \}$  defined in \eqref{e:CUSUM}, we are left to choose a threshold.    Denote for each $\kappa\ge0$, 
\begin{align}  
\thresh^*(\kappa) &= \argmin_{\thresh> 0} \{\kappa \MDE(\thresh) +  \MDD(\thresh)  \} 
\label{e:optThreshold}
\\
			J^*(\kappa) &= \min_{\thresh> 0} \{ \kappa \MDE(\thresh) +  \MDD(\thresh)  \}   
\label{e:optCost}
\end{align}
	We write CUSUM* to denote the CUSUM algorithm using the optimal threshold $\thresh^*(\kappa)  $. 
	
	\label{e:optThresholdCost}%
\end{subequations}%

In the following we survey results from \cite{coomey24b,coomey24c}
 only for the   independent   i.i.d.\ model.  
In addition, for the sake of brevity,  many explicit formulas below are shown for the conditional i.i.d.\ model with \emph{scalar} observations.
In this case it is assumed that the respective marginal densities $f_0$ and $f_1$ exist, so that the LLR becomes $L=\log(d\piX^1/d\piX^0)=\log(f_1/f_0)$, and 
			it is assumed that
$L$ is finite valued and integrable with respect to 
	both
	 $f_1$ and $f_0$.	See    \cite[Prop.~2.3]{coomey24b} for a proof of \Cref{t:barcApprox}.

\begin{subequations}
		
\begin{proposition}
\label[proposition]{t:barcApprox}
Suppose that (A1)-(A3) hold for the CUSUM test applied to the independent i.i.d.\ model. Then, the following approximations hold for the cost criterion \eqref{e:MDD+kappaMDE}:
\whamrm{(a)}   For each threshold we have the approximation $\barJ(\thresh, \kappa)  =   [1+o(1)]  \barJ_\infty( \kappa, \thresh) $\,, where
\begin{equation}
	\barJ_\infty( \kappa, \thresh)  =  \thresh  /m_1    +  \kappa \sqrt{\thresh}   \sqrt{2\pi  \gamma^2 }    \exp( -\thresh \thexp_+) 
	\label{e:barJinftySeptember}
	\end{equation}		
			in which  $ o(1) \to 0$  as $\thresh\to\infty$ and $\gamma^2 =  \Lambda_0''\, ( \thexp_+) /\thexp_+^3$.

\whamrm{(b)}  
			The optimal threshold $\thresh^*(\kappa)$ and cost $J^*(\kappa)$ admit the asymptotic approximations
\[
\thresh^*(\kappa) =  [1+o(1)]\,\barthresh_\infty^*( \kappa), 
\quad
			J^*(\kappa)  =  [1+o(1)]\,\barJ_\infty^*( \kappa) .
\]
in which  $ o(1) \to 0$  as $\kappa\to\infty$,  and
 \begin{equation}
	\barthresh_\infty^*( \kappa)  = \frac{1}{  \thexp_+ }          \log(\kappa)    
	\,,\quad 
	\barJ_\infty^*( \kappa) =
	\frac{1}{  m_1 }   \frac{1}{  \thexp_+ }        \log(\kappa) 
	\label{e:ApproxJH}
\end{equation}

\whamrm{(c)} 
Denote $F_r := F + r$ for any $r \in \Re$, resulting in approximate cost $\barJ_\infty^*( \kappa; F_r) $.   Its minimum over $r$ may be expressed  
\begin{equation}
\begin{aligned}
r^*& = \big(\expa - \Lambda_0(\thexp^*) \big)/\thexp^*
\\
&\textit{
where $\thexp^* = \argmin_\thexp [  \Lambda_0(\thexp)  -  \thexp \piX^1( F) ]$.}
\end{aligned}
\label{e:r-star}
\end{equation}

\whamrm{(d)} 
If the function $L$ satisfies  (A3), then
  $F^* = L + \expa$ minimizes  $\barJ_\infty^*( \kappa; F) $
among all functions $F$ satisfying (A1).    
\qed
\end{proposition}

\end{subequations}

The arguments in the proofs of (a) and (b) in \cite{coomey24b} begin with two observations:  1.~The cost of delay is easily approximated for this model:   After a change has occurred, the most likely path is linear with slope $m_1 >0$.   For a threshold $\thresh\gg 1$,  the delay  $ (\tstop -\tchange)_+$  is overwhelmingly likely to be close to $\thresh/m_1$.  2.~Approximation of the mean of $ (\tstop -\tchange)_-$ is based on  
	well-established large deviations theory for RRWs. Application of this theory rests on properties of the log moment generating functions introduced in (A3).	
		Parts (c) and (d) follow from \cite[Prop.~3.6]{coomey24b}.

\wham{Extension to more general statistics}

Assumption (A3) is valuable only for the independent i.i.d.\ model.   For models with memory, such as the independent Markov model, one must instead pose (A3) using the   \textit{cumulative} log moment generating functions \cite{coomey24b,coomey24c}.

Consider the independent Markovian model in which the two Markov chains have transition kernels $P_0, P_1$,  respectively,
and each  satisfies (DV3) (the Lyapunov function in 	\eqref{e:DV3body} may not be the same for each model).    The steady-state  marginal distributions are denoted $\piX^0,\piX^1$ respectively.    Extensions of the theory summarized in \Cref{t:barcApprox} tells us that  in application of CUSUM we should consider functions of the form $F_{n+1}=F(Y_n,Y_{n+1})$.    Assumption (A1) is then modified:  $m_0 = \Expect[F(X^0_n,X^0_{n+1})] <0 $  and  $m_1 = \Expect[F(X^1_n,X^1_{n+1})] >0 $, with expectations in steady-state. 

 Suppose there are transition densities $\{g_0,g_1\}$ with respect to a reference measure $\upmu$, so that $P_i(x,dz)=g_i(x,z)\,\upmu(dz)$,  and consider the transition log-likelihood ratio,
\begin{equation}
	L_\infty(x,z)=\log\Bigl(\frac{g_1(x,z)}{g_0(x,z)}\Bigr).
	\label{e:LLR_Markov}
\end{equation}
We find that this function satisfies the sign conventions required in (A1)  subject to   mild conditions on the Markov chains. These bounds are 
 justified using arguments similar  to 
\eqref{e:signcon1}.

A version of \Cref{t:barcApprox} holds, using $F^* = L_\infty + \expa$ in part (d)   \cite{coomey24b,coomey24c}.

\subsection{Q-learning in non-ergodic settings}
\label{s:Qregeb}

	In adopting   the CUSUM statistic as a SIS in an architecture for Q-learning, \Cref{t:Q} is not directly applicable since there is no underlying Markov chain that is ergodic:   
	Under (A1), \eqref{e:CUSUM} evolves as a transient random walk for $k\ge \tchange$.

We obtain an ergodic model through a small change in the objective \eqref{e:ObjPOMDPQCD}: 
It is assumed that $\Phi_0 \sim \muRegen$ for a probability measure $\muRegen$ supported on $\state_0$.   
Letting $\tregen$ denote the first return time to some set $\RegenState \subset \state$,
the modification of \eqref{e:ObjPOMDPQCD}  for a randomized policy $\feex $ is defined by    
\begin{equation}
J(\muRegen, \feex) = 
\Expect\Big[ \sum_{k=0}^{\tstop\wedge \tregen -1} \disc^k\, c_\circ(\Phi_k) 
			+ \disc^{\tstop\wedge \tregen}\, c_\bullet(\Phi_{\tstop\wedge \tregen}) \Big]
\label{e:POMDPregen}
\end{equation}
  \textit{This construction is introduced only for application of Q-learning}.   
It is justified assuming that data is generated using computer simulation.

The new objective permits the introduction of a sequence of regeneration times $\{  \RegenTime_k : k\ge 0 \}$ so that ergodicity can be assured.   
Set $\RegenTime_0=0$, and for $n\ge 0$ define $\RegenTime_{n+1}=\min\{k > \RegenTime_n : \Phi_k\in\RegenState \ \textit{or} \ U_k=1\}$.  
At each  $n\ge 1$
we take $\Phi_{\RegenTime_n}\sim\muRegen$ independently of the past. 	Consequently,  the sequence $\{\Phi_{\RegenTime_n}  : n\ge 1\}$ is i.i.d.\  with marginal $\muRegen$, and \eqref{e:POMDPregen} may be equivalently expressed,  for any $n\ge 0$,  
\[
J(\muRegen, \feex) = 
\Expect\Big[ \sum_{k=\RegenTime_n}^{\RegenTime_{n+1} -1} \disc^{k - \RegenTime_{n} }\, c_\circ(\Phi_k) 
			+ \disc^{\RegenTime_{n+1} -\RegenTime_{n} }\, c_\bullet(\Phi_{\RegenTime_{n+1} }) \Big]
\]

The  projected Bellman equation \eqref{e:pBe} is modified slightly through a modification of the temporal difference:
\begin{equation}
	\TD^\theta_{k+1} \eqdef -\,Q^\theta_k + c(\Phi_k,U_k)
		+ \disc\,(1-U_k)\,\ind_{\RegenState^c}\,\uQ^\theta(\InfoState_{k+1}) .
	\label{e:TDregen}
\end{equation}
where $\ind_{\RegenState^c} \eqdef \ind\{\Phi_k \notin \RegenState\}$. The associated Q-learning algorithm with discount factor $\disc<1$ is precisely \eqref{e:POMDPQ} using \eqref{e:TDregen}.
	Extension of \Cref{t:Qgamma1a} for the associated Q-learning algorithm with $\disc=1$ is straightforward.

	The required assumptions for convergence are in general weaker than those in \Cref{t:Qgamma1a},  based on the rank of the sum of two matrices: instead of requiring $\Sigma_{k+1\mid k}$ to be full rank, it suffices that $\Sigma_{k+1\mid k}+ M^\RegenState$ be full rank, where
\[
\begin{aligned}
	\Sigma_{k+1\mid k} &= \Expect_\piPhi[  ( \psi_{(k+1)} - \psi_{(k+1 \mid k)}  ) ( \psi_{(k+1)} - \psi_{(k+1 \mid k)}  ) ^\transpose]
	\\
		M^\RegenState &=  \Expect_\piPhi[  \ind_{\RegenState^c} (\Phi_k) \psi_{(k+1 \mid k)} \,\psi_{(k+1 \mid k)}^\transpose]
\end{aligned}
\]
The first was previously defined in \eqref{e:condCov}, with
	$\psi_{(k+1 \mid k)} =  \Expect[\psi(\InfoState_{k+1}, U_{k+1})\mid \Phi_k]$.

The proof of the following is essentially identical to   the proof of \Cref{t:Qgamma1a}.
\begin{proposition}
	\label[proposition]{t:Qgamma1}   
		Consider the mean flow vector field \eqref{e:pBe} using the temporal difference sequence 
	\eqref{e:TDregen}:
		$\barf(\theta)  =  $
	\[
	\Expect\big[  \psi_{(k)} \{- \psi_{(k)}^\transpose \theta   +  c(\Phi_k, U_k)  
		+  \disc (1-U_k)   \ind_{\RegenState^c} (\Phi_k)  
	\uQ^\theta(\InfoState_{k+1} )   \}  \big]
	\]
		Suppose the assumptions of \Cref{t:Q} hold.  Then, the conclusions of  \Cref{t:Q} hold if  in addition  $\Sigma_{k+1\mid k} + M^\RegenState$  is full rank.    
		\qed
\end{proposition}

\section{Numerical Results} 
\label{s:Qqcd}

This section concerns design and evaluation of   Q-learning algorithms in the regenerative setting of \Cref{s:Qregeb}, so in particular we are interested in  the performance criterion \eqref{e:POMDPregen}.   
 We begin by describing elements common  to algorithm design for these examples.  
 
For the (possibly multidimensional) SIS 
process $\{\InfoState_k\}$ we fix a large constant $\upeta$  to define 
 $\RegenTime_{n+1}=\min\{k > \RegenTime_n :  \|\InfoState_k\|\ge\upeta   \ \textit{or} \ U_k=1\}$. 
The input used for training in the Q-learning algorithms was chosen i.i.d., with  $p_i=\Prob_\piPhi\{ U_k = i \}  =1/2$ for $i=0,1$.   
 Further elements are summarized in the following:

\notes{
I don't think we need this here, but let's see:
\\
 Results are surveyed for both i.i.d.\ and Markovian independent models as  defined in  \Cref{s:POMDPQCD}. 
}

 \wham{Acceleration}

A matrix gain algorithm was implemented, known as \textit{Zap Q-learning} \cite{chedevbusmey20a,devmey17b}, 
\begin{equation}
	\theta_{n+1} =  \theta_n + \alpha_{n+1}  G_{n} \zeta_n  \TD_{n+1}
	\label{e:POMDPQ_Zap}
\end{equation}
	with  $G_{n} $   an approximation of  $- [\barA(\theta_n)]^{-1}$,  where 
\[
\barA(\theta) \eqdef \partial_\theta \barf\, (\theta) = -R + \disc  p_0  \Expect_\piPhi[ \psi^0(\InfoState_k) \psi(\InfoState_{k+1}, \fee^\theta(\InfoState_{k+1}) ) ^\transpose ]
\]
See \cite[Section~8.3]{CSRL} for details on the recursion defining $\{G_n \}$.  The convergence rate of Q-learning was observed to be extremely slow using the scalar-gain algorithm,  even with the application of averaging \eqref{e:PR}, 
whereas Zap Q-learning performed far better;    explanation is provided in   \cite{CSRL,chedevbusmey20b}.     
	
Successful implementation requires invertibility of $\barA(\theta)    $  for ``most''  $\theta\in\Re^d$  and the step-size choice   $\alpha_n = \min(\alpha_0,1/n)$ for any $\alpha_0>0$ \cite{chedevbusmey20b,chedevbusmey20a}. A proof of invertibility for all $\theta$ follows as in the proof of stability of the mean flow
(recall \Cref{t:QstableMeanFlow}).

For any algorithm, the   \textit{asymptotic covariance} is defined by
\begin{equation}
\SigmaTh =
		\lim_{n\to\infty} n \Expect[  \, \tiltheta_n \{\tiltheta_n  \}^\transpose]
\label{e:avarPR}
\end{equation}
where 	$\tiltheta_n \eqdef \theta_n-\theta^*$.    The asymptotic covariance obtained using Zap Q-learning is precisely the optimal $\SigmaPR $ appearing in \eqref{e:PR}.    The conclusion is a minor extension of  the main result from \cite{borchedevkonmey25}, which also asserts that the
Central Limit Theorem (CLT) holds:  $Z_n \eqdef \sqrt{n} \tiltheta_n  \darrow N(0, \SigmaPR)$, where the convergence is in distribution.

The optimal asymptotic covariance admits the representation  
\begin{equation}
\SigmaPR = G^* \Sigma_\Delta [G^*]^\transpose 
\label{e:ChenPolyakRupertCov}
\end{equation}
  in which $G^* = - [\barA(\theta^*) ]^{-1}   $  and $\Sigma_\Delta $ is the power spectral density matrix
$
  \Sigma_\Delta  =  \sum_{k=-\infty}^\infty \Expect[\Delta_0\Delta_k^\transpose]
 $,
where $\Delta_k =  f(\theta^*, \qsaprobe_{k+1})$ and the expectations are in steady-state.   In practice we cannot easily estimate     $  \Sigma_\Delta$, but CLT theory justifies estimation of the asymptotic covariance through alternative means.

 \wham{Considerations for one dimensional SIS.}
When using a one-dimensional RRW as a SIS in Q-learning we 
define a threshold as follows:     
\begin{equation}
\clH(\theta) \eqdef \inf \{s  :  Q^\theta(s,0) > Q^\theta(s,1) \}
\label{e:clHtheta}
\end{equation}  
In most cases we find the policy  \eqref{e:fee_theta_HTorQCD} coincides with the threshold policy, in the sense that   $\fee^\theta(s)=\ind\{ s \ge \clH(\theta) \}$.     To compare the resulting threshold with the optimal we performed long run-lengths to collect estimates of the crucial data
\begin{equation}
\thresh^*(\kappa)   \,,   \ \   \clJ(\thresh;\kappa)  \,, \quad  \kappa>0\, \ \thresh>0
\label{e:FixedStats}
\end{equation}
with $\thresh^*(\kappa) $ defined in \eqref{e:optThreshold},   and $\clJ $  the approximation of \eqref{e:MDD+kappaMDE} defined in \eqref{e:POMDPregen}
(see \Cref{s:appndx_exp} for details).  
While only estimates, we will use the notation $\thresh^*$, $  \clJ$ without further comment in the discussion that follows.

 \wham{Evaluation}
In applying Q-learning within the regenerative setting of \Cref{s:Qregeb}, for given $N\ge 1$ we   let $\EpsLength=\EpsLength(N)\ge N$ denote the total number of samples $(\InfoState_k, U_k)$, so that  $\hatheta = \theta_{\EpsLength}$ is the final estimate, which   defines a policy $ \hafee \eqdef \fee^{\hatheta}$ via \eqref{e:fee_theta_HTorQCD}.

The performance of an algorithm is multi-faceted.   Of course,  we want $\fee^{\theta^*}$ to be nearly optimal, with  $\theta^*$ the limit of the algorithm.   
A good algorithm will also rapidly obtain a good estimate of $\theta^*$.  
In the absence of tight non-asymptotic  error bounds, we adopt as a  measure of reliability 
the asymptotic covariance \eqref{e:avarPR}, which may be estimated  using the   \textit{batch means method}.   This requires  $M$ independent runs, resulting in  $\{ \theta_{\EpsLength^i} ,  \EpsLength^i: 1\le i\le M\}$,
with the average of the $M$ parameter estimates   denoted $\bartheta $.    
The empirical covariance of $\{ Z^i = \sqrt{\EpsLength^i} [ \theta_{\EpsLength^i} - \bartheta ] : 1\le i\le M \}$ provides the batch means 
 estimate of $\SigmaTh$.     
 Consistency (as $M\to\infty$ and then $ N\to\infty$) follows from \eqref{e:avarPR}  and the Law of Large Numbers for i.i.d.\ random vectors.   
 
\begin{figure}[h]
	\centering		 
	\includegraphics[width=0.7\hsize]{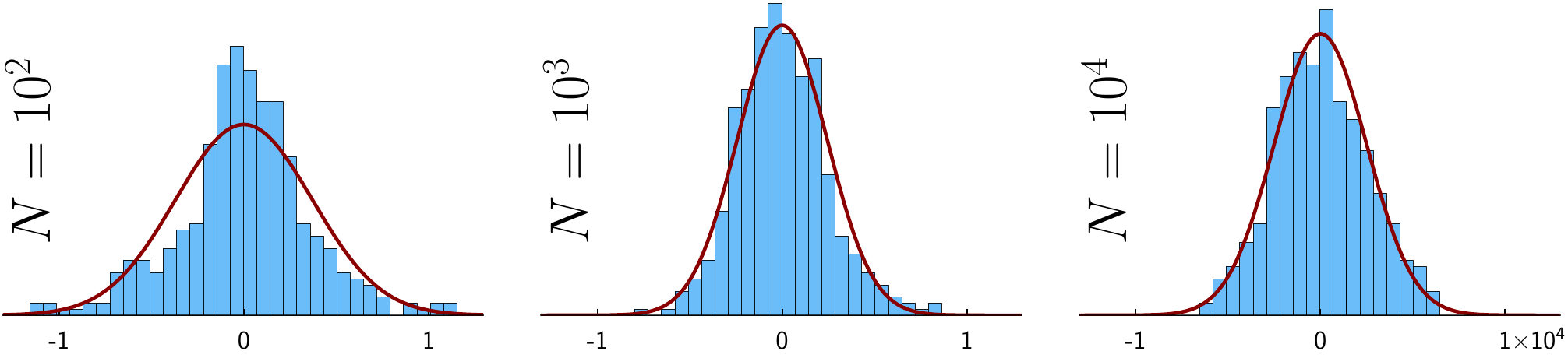}
	\caption{Histograms of $\{Z^i_1 : 1\le i\le M\}$  for three values of $N$. }
	\label[figure]{f:Hists}
\end{figure}

 An example of the CLT is provided in  \Cref{f:Hists} for \textbf{QCD Model 1}  (details may be found below).  It is evident that the Gaussian approximations for the two larger values of $N$ are highly consistent.

Models 1 and 2 below concern a one-dimensional RRW as a SIS, in which case it is of  greater interest to approximate the variance of $\clH(\theta_{\EpsLength} ) $.   Under the assumption that $\clH(\theta) $ is twice continuously differentiable in a neighborhood of $\theta^*$, 
there are two possibilities: 

 \wham{1.  Ideal:   $\clH(\theta^*) = \thresh^*(\kappa)  $.}  In this case,  a second order Taylor series combined with  \eqref{e:avarPR}  implies that $  \Expect[  ( \clH( \theta_{\EpsLength} ) -\thresh^*(\kappa)  )^2]$ vanishes as $O(1/N^2)$ as $N\to\infty$, which is far faster than anticipated by the CLT.     
 
 Approximation of the variance of the threshold can be obtained via the CLT for the parameter estimates.
 On denoting  $\{ \clZ^i =  {\EpsLength^i} [ \clH( \theta_{\EpsLength^i} ) - \bartheta ] : 1\le i\le M \}$, 
 we have    for each $i$ the distributional limit $\clZ^i \darrow \clZ_\infty $  as $N\to\infty$, where   $\clZ_\infty =  W M_h W $ in which $M_h = \nabla^2\clH(\theta^*) $  and $W\sim N(0, \SigmaTh)$  (the random variable  $\clZ_\infty $ has a \textit{generalised chi-squared distribution}).   

 \wham{2. $\clH(\theta^*) \neq \thresh^*(\kappa)  $.}  A   first-order Taylor series implies a CLT for 
$\{ Z^i =  \sqrt{\EpsLength^i} [ \clH( \theta_{\EpsLength^i} ) - \bartheta ] : 1\le i\le M \}$,   with limiting distribution $N(0, \sigma^2_h)$ as $N\to\infty$,  where
$\sigma^2_h =   [\nabla\clH\, (\theta^*) ]^\transpose \SigmaTh  \nabla\clH\, (\theta^*) $.

\smallskip

\begin{figure}[h]
	\centering	 
	\includegraphics[width=.4\hsize]{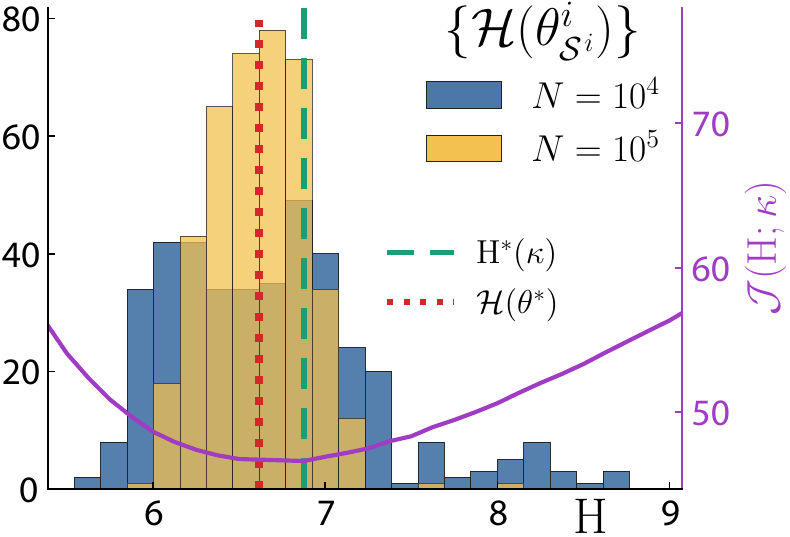}
	\caption{Histograms of policy thresholds over independent experiments for two values of $N$. Also shown are the minimizing threshold for CUSUM* and the threshold from the final Q-learning run using a larger $N$.  }	
	\label[figure]{f:CostHistos}
\end{figure}

An illustration based on \textbf{QCD Model 1a}  is contained in
\Cref{f:CostHistos}.  The two dashed lines indicate that $\clH(\theta^*) \neq \thresh^*(\kappa)  $,  but from the plot of $\clJ$ it is seen that 
the threshold $\clH(\theta^*) $ results in a policy that is very nearly optimal  (recall the discussion surrounding \eqref{e:FixedStats} for notation).   
Also shown are the un-normalized histograms of the thresholds  $\{\clH(\theta_{\EpsLength^i}^i) : 1 \le i \le M\}$  for two values of $N$.   
The Taylor series approximation is not very useful for $N = 10^4$: the histogram is far from Gaussian (two outliers with $\clH(\theta_{\EpsLength^i}^i)>9$ were removed).   
 
 A Gaussian approximation appears justified for $N=10^5$, though the high volatility suggests that $\|  \nabla\clH\, (\theta^*) \|$ is large.

\wham{Basis selection}  

The basis for the 
	linear 
function class ${Q^\theta =\theta^\transpose \psi}$ 
was chosen to respect the assumptions of \Cref{t:Q}.  In particular, under the separability assumption,  
we may impose without further loss of generality the form  
$\psi(s,u) = [ (1-u) \psi^{(0)}(s) ; u \psi^{(1)}(s) ]$.
For each $j=0,1$ and given $d\ge 2$,  the components of  $\psi$ were defined using Gaussian radial basis functions (RBFs),
\begin{align}
	\psi^{(j)}_i(s) = \exp\left( - \half \sigma_i^2 \| s - \mu_i \|^2    \right)\,, \qquad 1 \le i \le K =d/2 \, ,
	\label{e:basis2}
\end{align} 
where     the
centers $\{ \mu_i \}_{i=1}^K$ were obtained through $k$-means clustering applied to a collection of CUSUM samples $\{s^{(m)} \}_{m=1}^M$  from the QCD model considered in preliminary Q-learning experiments.    
 For each $\mu_i$ define the width parameter  $\sigma_i = b  \min_{l \neq i} \|\mu_i - \mu_l\|$
where scalar $b > 0$ is a design choice.

Parameter values used for \eqref{e:basis2} are summarized in \Cref{s:appndx_exp}.  

We also tested bases based on binning,  of the form   $\psi_i(s,u) = \ind\{s\in S_{k_i} ,   u = u^i\} $  for a collection of intervals $ \{S_j \}$ and input values $\{u^j\} $.     The results were often poor, as might be anticipated by examining the structure of the policy defined in \eqref{e:fee_theta_HTorQCD}.    In the one-dimensional SIS examples considered, the policy  $\fee^\theta$   typically falls in the class of threshold policies, whose threshold  $\clH(\theta)$ defined in \eqref{e:clHtheta} must be  a boundary of one of the intervals $ \{S_j \}$.    This  in general implies we require a large number of bins to approximate the best policy.

\begin{figure}[h]
	\centering		
	\includegraphics[width=0.52\hsize]{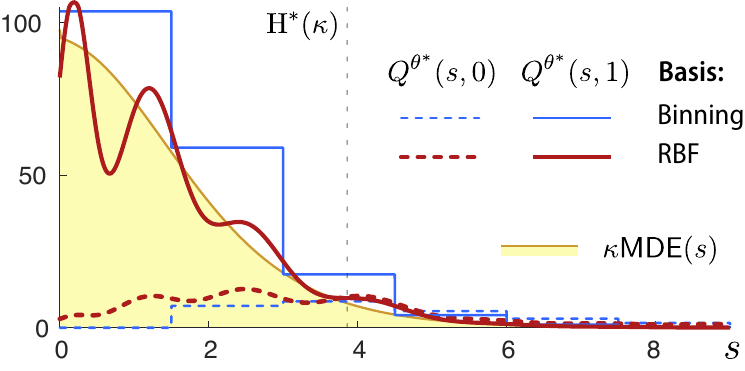}
	\caption{Q-function approximations obtained in two experiments, differentiated by the choice of basis:  1.~binning,  and  2.~RBF defined in    	\eqref{e:basis2}}
	\label[figure]{f:Qsolutions}
\end{figure}

 \Cref{f:Qsolutions} compares results obtained using binning and using a RBF  for \textbf{QCD Model 1a} for $\kappa=2$.  
 Based on the definition \eqref{e:fee_theta_HTorQCD},  from these plots we conclude that   $\hafee \eqdef \fee^{\theta^*}$ is a threshold policy for either choice of basis.  
For the choice of RBF we found in most cases that the threshold approximated the optimal threshold $\thresh^*(\kappa)  $ that defines CUSUM*
  (recall 	\eqref{e:optThresholdCost}
for notation).     See  \Cref{s:QCDmod1} for explanation of the plot of 
 $\kappa \MDE(s)$   shown in the figure. 

The remainder of this section concerns Q-learning for three  instances of the independent model described in \Cref{s:POMDPQCD},
 distinguished by the assumptions on the observations and the change time, and on the features used for Q-learning.   
Results are surveyed only for experiments performed using RBFs \eqref{e:basis2} to define $\psi$.

\subsection{QCD Model 1.   Conditional i.i.d.\ model} 
\label{s:QCDmod1}

\wham{Statistical model}

The two stochastic processes in \eqref{e:QCDmodel}  are taken i.i.d.\ Gaussian with
$\preObs_k \sim N(0,\sigma^2) = \preDens$ and 
$\postObs_k \sim N(\mu_1,\sigma^2)  = \postDens$ for each $k$,  with $\mu_1 = 0.5$ and $\sigma = 1$.    	Two cases for the change time $\tchange$ are considered, each satisfying $\expa=0.02$ in  \eqref{e:hazardAss}:

\whamit{Shiryaev's model: } The  distribution of $\tchange$ is geometric with parameter $\expa$.

\whamit{Mixed change time: } The distribution is a mixture:  $\tchange \sim \text{geo}(\expa)$ with probability $0.25$; else $\tchange \sim \text{geo}(0.2)$.

Shiryaev's model is a convenient starting point since the optimal test is available in a very simple form.

\begin{wrapfigure}[10]{r}{0.25\hsize}
	\vspace{-0.0em}					 %  <--  Playing with spacing around wrapfig (!)
	\centering     
	\includegraphics[width=\hsize]{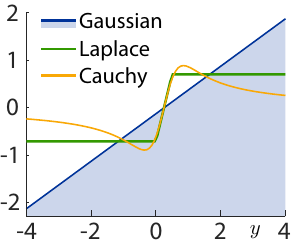}
	\vspace*{-5.5mm} % moving caption closer to fig
	\caption{$\surL$ for Models \modelAone-\modelAthree.}
	\label[figure]{f:3LLRs}
	\smallskip
	
\end{wrapfigure}
	
	\wham{CUSUM architecture}

Three functions were tested to define  the CUSUM statistic,
each of the form  $F=\surL+r^*$ with different (possibly) ``mismatched'' densities for $\surL=\log(\surf_1/\surf_0)$,
and $r^*$ defined in \Cref{t:barcApprox}~(c).
	  Plots of the three LLRs defined below are shown in     \Cref{f:3LLRs}. 
\\
\textbf{\modelAone: Gaussian.} $\surf_0$ is $N(0,1)$ and $\surf_1 $ is $N(\mu_1,1)$ with $\mu_1 = 0.5$. Hence $\surL=L$.
\\ 
\textbf{\modelAtwo: Laplace.}   $\surf_0$ is Laplace$(0,b)$ and $\surf_1 $ Laplace$(\mu_1,b)$ with $\mu_1 = 0.5$ and $b = \sqrt{\sigma^2/2}$ (matching second order statistics).
\\
\textbf{\modelAthree: Cauchy.} $\surf_0$ is Cauchy$(0,\CauchyScale)$ and $\surf_1 $ Cauchy$(x_1,\CauchyScale)$ with $x_1 = 0.5$ and $\CauchyScale$ chosen so that the Gaussian and Cauchy CDFs evaluated at $\sigma = 1$ are equal.

To obtain   $r^*$  we must solve \eqref{e:r-star}.    In \textbf{Model \modelAone}, an application of \Cref{t:barcApprox}~(d) gives $r^* = \expa$.  
For mismatched \textbf{Models \modelAtwo} and \textbf{\modelAthree}, the value of $r^*$ was computed as 0.031 and 0.036 by an application of \Cref{t:barcApprox}~(c).

\begin{subequations}

\wham{Ideal performance and approximations}

\Cref{f:CostAndThresholdApproximations} depicts CUSUM* values $\{ \thresh^*(\kappa),  J^*(\kappa)  \}$ 
defined in \eqref{e:optThresholdCost} for this model.   We find that the approximations found in \Cref{t:barcApprox} are poor in this example for this range of $\kappa$.   However, consider  the shifted values
\begin{align}
\barthresh_{\tiny{\sf s}}~(\kappa)  &=
\barthresh_\infty^*(\kappa) - \barthresh_\infty^*(100) + \thresh^*(100)
\label{s:shiftH}
\\
\barJ_{\tiny{\sf s}}(\kappa)  &=  \barJ_\infty^*(\kappa)  -
\barJ_\infty^*(100)    +  J^*(100)
\label{s:shiftJ}
\end{align}
	This ensures that the approximations \eqref{e:ApproxJH} coincide with the estimates of  \eqref{e:optThresholdCost} at $\kappa=100$.

	\label{s:shiftHJ}%
\end{subequations}
	
These values are also shown in  \Cref{f:CostAndThresholdApproximations}, and closely approximate the CUSUM* values.
Hence the error $ \barthresh_\infty^*(\kappa) - \thresh^*(\kappa) $ is nearly constant over the entire range.  

	Unfortunately, the constant value is large, which is one motivation for learning techniques to obtain a near-optimal policy even in a   one-dimensional SIS setting.

\smallskip

We turn next to results obtained using Q-learning.

\smallskip

As remarked earlier, the policy  obtained from Q-learning resulted in a threshold policy in almost all cases.    Therefore, the best performance from Q-learning can be no better than CUSUM*,  which uses the optimal threshold \eqref{e:optThreshold}.

\begin{figure}[h]
	\centering	 
	\includegraphics[width=0.755\hsize]{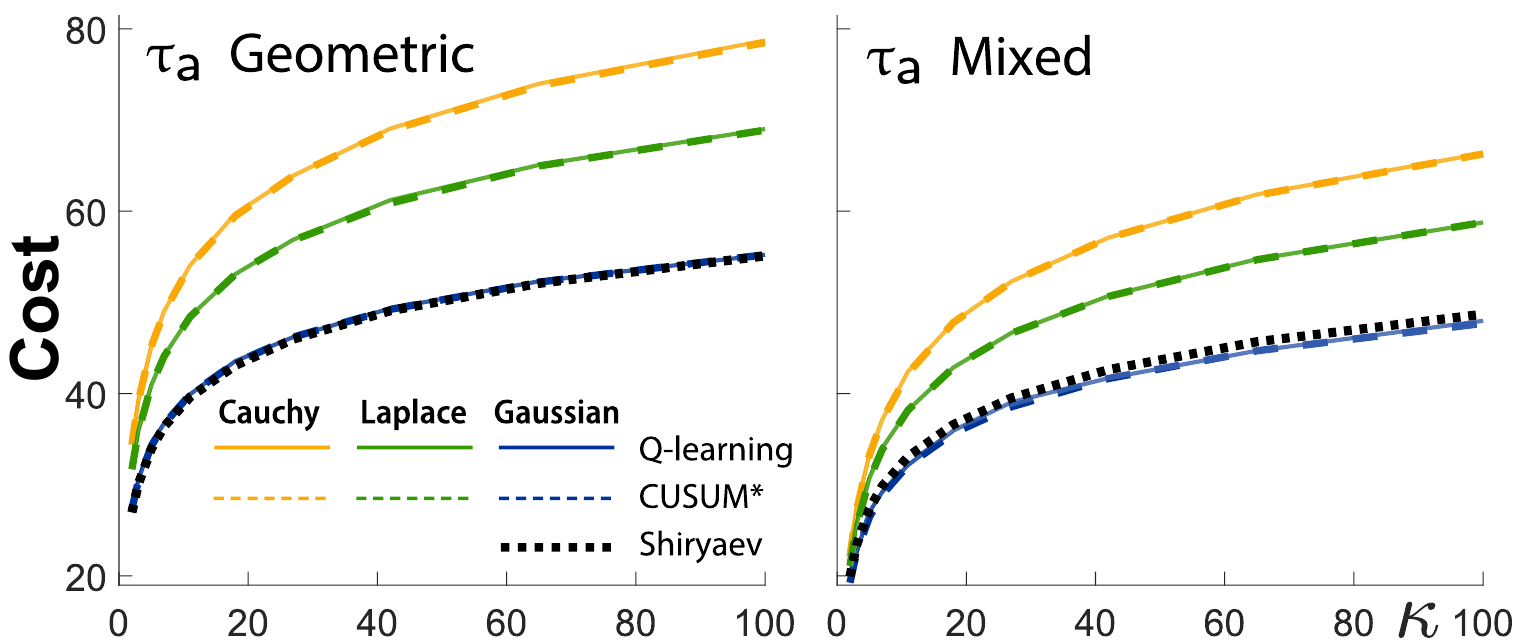}
	\caption{Q-learning performance for QCD Model 1 considering two cases: Shiryaev's model on the left and mixed change time on the right. }		
	\label[figure]{f:PerformanceAll3}
\end{figure}

\whamit{Performance evaluation.} The left-hand side of \Cref{f:PerformanceAll3} shows the average cost obtained from the Q-learning policy $\fee^{\theta^*}$ for   \textbf{\modelAone}. Also included are CUSUM* and the \textit{optimal test}---Shiryaev's defined in \Cref{s:POMDPQCD}.   The average cost curve obtained through Q-learning matched remarkably close to CUSUM* for all $\kappa$ tested.

Experiments were repeated for \textbf{Models \modelAtwo} and \textbf{\modelAthree}. The simulation environment for Q-learning remains the same for each case, except that the SIS  $\{\InfoState_k\} $ differs based on the respective LLRs, plotted in \Cref{f:3LLRs}.

Revisiting \Cref{f:CostAndThresholdApproximations}, we bring attention to the  poor performance of CUSUM* for mismatched \textbf{\modelAtwo} and \textbf{\modelAthree} compared to \textbf{\modelAone}. In this sense, we noted consistent agreement between CUSUM* and results from Q-learning:
In all three cases, the policy yielded  average cost within 3\% of CUSUM*. The right-hand side of \Cref{f:PerformanceAll3} shows similar performance outside of Shiryaev's model. Notably Q-learning outperforms Shiryaev's test, which is no longer optimal for the mixed change time distribution.

\whamit{Illustration of  $L_2$ theory}
Recall that \Cref{t:Q}~(b) establishes that $Q^{\theta^{*}}(s,1) /\kappa $    
defines the projection of the mean detection eagerness onto the random variables spanned by basis elements $\{\psi^{(1)}_i \}$.  
We let
$\MDE(s)$ denote the \textit{conditional mean detection eagerness},  conditioned on $\InfoState_k = s$.    
That is,      
\begin{align*}
	\tfrac{1}{\kappa} Q^{\theta^{*}}(s,1)  & =      \widehat{\Expect}\big[ (\tstop -\tchange)_-    \mid \psi^{(1)}_{(k)} = \psi^{(1)}(s)  \big] 
	\\
	\MDE(s)  &=  \Expect\big[(\tstop -\tchange)_-  \mid \InfoState_k = s \big]  
\end{align*}
where $\widehat{\Expect}$ denotes projection onto the span of $\{\psi^{(1)}_i \}$ within the Hilbert space $L_2(\piPhi)$.

Shown in \Cref{f:Qsolutions}  is a   plot of $\kappa \MDE (s)$ as a function of $s$  (obtained via Monte-Carlo) for \textbf{Model \modelAone}.  
The approximation  $Q^{\theta^{*}}(s,1) \approx \kappa \MDE(s)$ is far more accurate with the RBF basis \eqref{e:basis2} when compared to binning.

\whamit{Convergence rates.}

We return to the histograms shown in \Cref{f:Hists,f:CostHistos} and the surrounding discussion. Each histogram was generated from $M=400$ independent experiments using Shiryaev's model and $\kappa=27$. 
Only the first component of the histogram is shown for \Cref{f:Hists}, giving an estimate of $\SigmaTh(1,1)$.
What is crucial here is that the estimate of this value 
is nearly identical for the three values of $N$ chosen.    Similar results were observed for estimates of the other diagonal entries of $\SigmaTh$.
In this experiment we find that
$N = 10^3$ provides a reasonable estimate of the variance of $Z^i = \sqrt{\EpsLength^i} [ \theta_{\EpsLength^i} - \bartheta ]$ for each $i$.

\wham{Variations of the statistical model.}

Variations on \textbf{Model 1} were tested in which the two processes in \eqref{e:QCDmodel} were taken non-Gaussian,
maintaining the  three  CUSUM architectures of \textbf{Model 1}.   Selected findings are summarized below:

\whamb   $\preObs_k \sim \text{Laplace}(0,b)$ and $\postObs_k \sim \text{Laplace}(\mu_1,b)$ so that \textbf{\modelAtwo} matches the observation densities in the sense of \Cref{t:barcApprox}~(d).
For all $\kappa$ tested, Q-learning with SIS defined by   \textit{either}  
  \textbf{\modelAtwo} \textit{or} \textbf{\modelAthree} produced average cost close to CUSUM*, while the SIS based on \textbf{\modelAone} showed comparable performance only for $\kappa \le 18$.

\whamb  $\preObs_k \sim \text{Cauchy}(0,\CauchyScale)$ and $\postObs_k \sim \text{Cauchy}(x_1,\CauchyScale)$ so that \textbf{\modelAthree} matches the observation densities
in the sense of \Cref{t:barcApprox}~(d).  
 In this setting \textbf{\modelAone} violates assumption (A1) in \Cref{s:asy} and was not used.  Q-learning with SIS based 
 on 
\textbf{\modelAtwo} or \textbf{\modelAthree}
 again
   yielded average cost close to CUSUM* for all $\kappa$ tested.

The fact that the function $F=\surL + r^*$ obtained from either 
  \textbf{\modelAtwo} or \textbf{\modelAthree}
  gave similar performance is not surprising given the plots of 
    $\surL$ shown in \Cref{f:3LLRs}.
	
\subsection{QCD Model 2.   Conditional Markov model} 
\label{s:model2}

\wham{Statistical model}

The two processes are defined by a scalar linear Markov model $\Phi^i_{k+1}=A^i \Phi^i_k + W^i_{k+1}   $, subject to the following:

\begin{wrapfigure}[13]{r}{0.375\hsize}
	\vspace{-1.02em}					 %  <--  Playing with spacing around wrapfig (!)
	\centering     
	\includegraphics[width=\hsize]{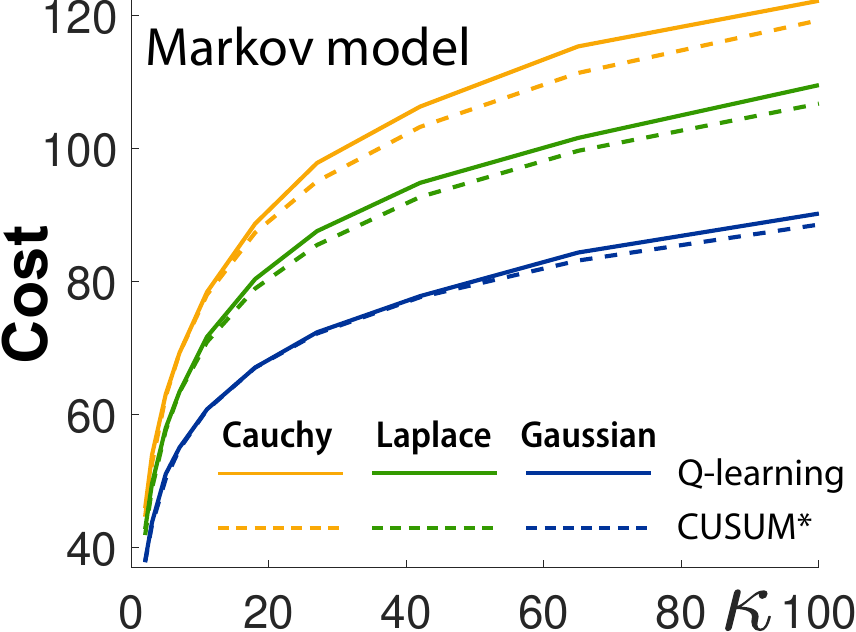}
	\vspace*{-5.5mm} % moving caption closer to fig
	\caption{Q-learning performance for   Model 2.}
	\label[figure]{f:PerformanceMarkovAll3}
\end{wrapfigure}

\whamc The change time $\tchange$ has a geometric distribution with parameter $\expa=0.02$.

\whamc $\{W^i_k : k\ge 1\}$ is i.i.d.\ $N(0,1)$ for each $i$.  
Consequently, only the dynamics change at time $\tchange$.

\whamc $A^0=0.8$ and $A^1=0.5$.

Under these assumptions we obtain from \eqref{e:LLR_Markov} the expression $L_\infty(x,z)=a x^2 - b xz$ with $a=([A^0]^2-[A^1]^2)/2=0.195$ and $b=A^0-A^1=0.3$.

\wham{CUSUM architecture}
Based on the theory surveyed in \Cref{s:asy} we chose 
$F_{n+1}=F(Y_n,Y_{n+1})$
in \eqref{e:SR},  with $F = \surL + r^* = \log(\surg_1/\surg_0) + r^*$.    
Three choices of $\surL$ were tested, each of the form  $\surg_i(x,z) = \surf(z - A^i x )   $ for each $i$,
with  scalar   $r^*$ obtained according to \Cref{t:barcApprox}~(c):  
\\
\textbf{\modelBone:   Gaussian.}   
$\surf$ is the $N(0,1) $ Gaussian density, so that       $r^* = \expa$. 
\\ 
\textbf{\modelBtwo: Laplace.}   $\surf$ is the Laplace$(0,b)$ density with $b = \sqrt{1/2}$, giving   $r^*=0.031$.  
\\
\textbf{\modelBthree:  Cauchy.} $\surf$ is the Cauchy$(0,\CauchyScale)$ density with $\CauchyScale  $ defined as in \textbf{Model~\modelAthree} and   $r^*=0.036$.

\whamit{Performance evaluation.} 
The plots shown in \Cref{f:PerformanceAll3} illustrate that the performance obtained from Q-learning is nearly optimal (though the performance gap is larger than observed in  \textbf{QCD Model 1}).

\whamit{Convergence rates.}

The histograms for $Z^i = \sqrt{\EpsLength^i} [ \theta_{\EpsLength^i} - \bartheta ]$ closely resembled those shown in \Cref{f:Hists}, but with standard deviations approximately doubled. Despite the higher variability in the parameter estimates, histograms of the threshold values $\{\clH(\theta_{\EpsLength^i}^i)\}$ showed much lower variance when compared to results from \textbf{Model 1} in \Cref{f:CostHistos}.

\subsection{QCD Model 3.   Multidimensional features} 

\notes{AC: Edited Model 3 lead-in for clarity. Could be better. Wasn't sure which sentences you wanted flipped, though.
\\
I've lost track,  but keep this note to remind me in a future revision.    We will define "notes" as null before submission.}

Results so far show that Q-learning performs well compared to CUSUM*, but the resulting cost may be large depending on the choice of SIS. Moreover, construction of $F^*$ requires substantial knowledge of the statistics for $\bfpreObs$, $\bfpostObs$ in \eqref{e:QCDmodel}.  This motivates features of multiple RRWs to define a multidimensional SIS.   
	
\wham{Statistical model}   The  observation model was \textbf{QCD Model 1} with geometric change time.

	\wham{CUSUM architecture} The SIS is defined by a pair of functions for CUSUM,  giving $\InfoState_n \in\Re^2$
	 for each $n$.

\wham{\modelCone: Laplace and Cauchy.}  Pair $(F^b,F^c)$ is defined by \textbf{Models  \modelAtwo} and \textbf{\modelAthree}.
\wham{\modelCtwo: Gaussian pair.} $(F^*,F^d)$ with $F^*$ defined by \textbf{Model  \modelAone}. A mismatched  Gaussian LLR  supposing $\surf_1$ is $N(0.1,1.4)$ defines $F^d$.
	
\begin{figure}[h]
	\centering	 
	\includegraphics[width=.75\hsize]{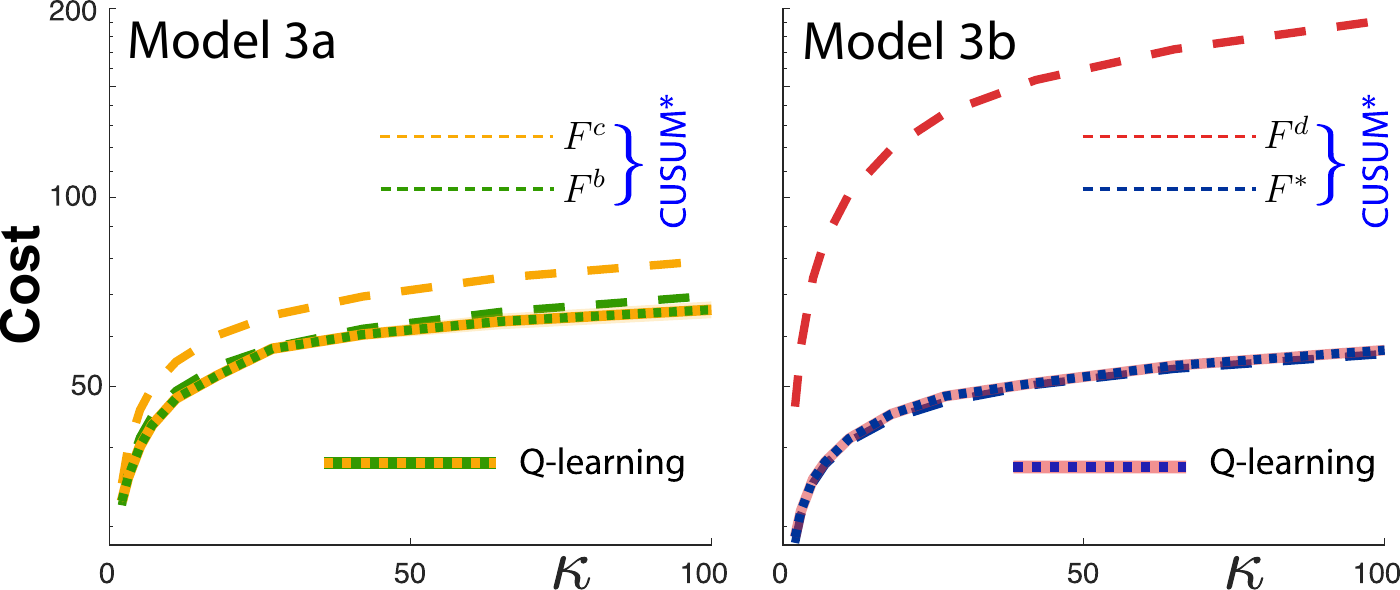}
	\caption{Q-learning performance for QCD Model 3 compared to CUSUM* for each SIS component. }		
	\label[figure]{f:PerformanceQ2D}
\end{figure}

	We constructed $F^d$ in \textbf{Model \modelCtwo}  to give especially poor performance compared to $F^*$, as illustrated on the right-hand side of \Cref{f:PerformanceQ2D}.

\begin{figure}[h]
	\centering	 
	\includegraphics[width=.9\hsize]{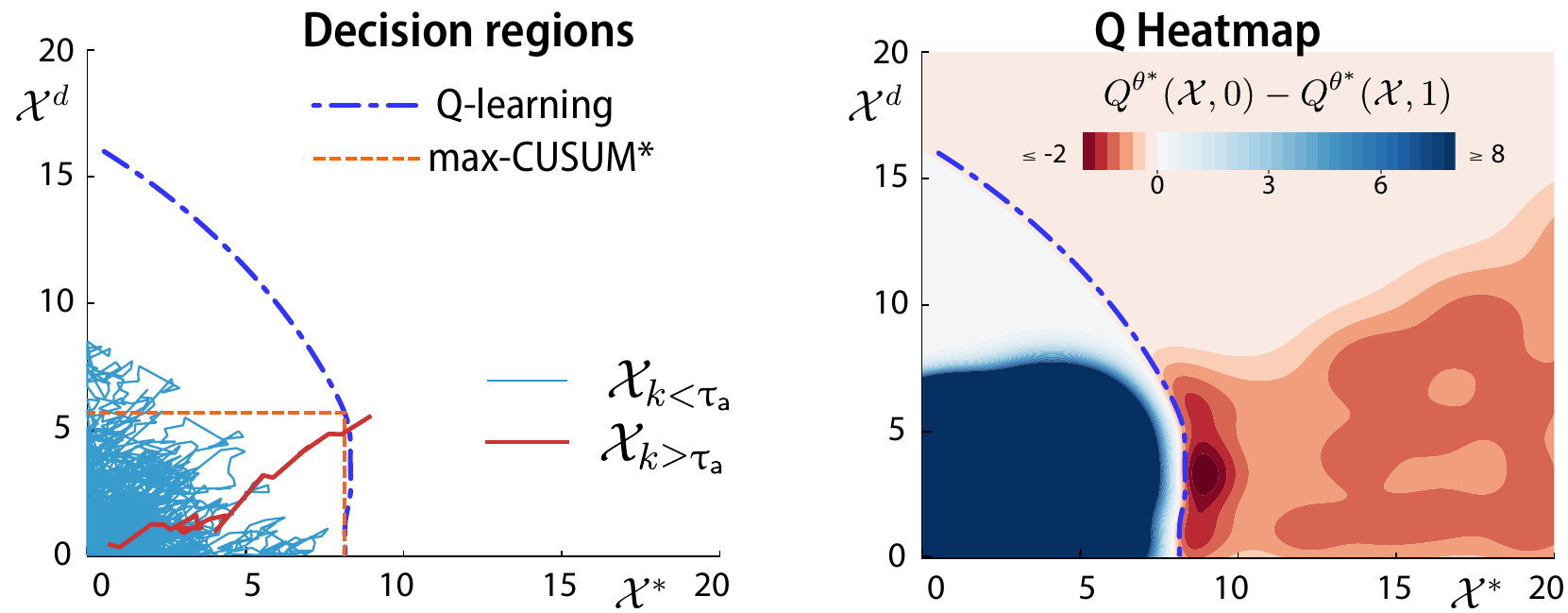}
	\caption{The left hand side shows the decision regions obtained from both Q-learning and max-CUSUM*,  along with a typical sample path. The right hand side depicts a policy heat map induced by the Q-function.}		
	\label[figure]{f:Q2D_region+trajs}
\end{figure}

Let \textit{max-CUSUM*}   denote the decision region defined by a box constraint,  with boundaries determined by CUSUM* thresholds.   An example may be found in
  \Cref{f:Q2D_region+trajs}, in which the box defines the acceptance region for the stopping rule.  The figure is based on  \textbf{\modelCtwo} with $\kappa=100$, for which the max-CUSUM* decision region is defined by $(\thresh^*(100), \thresh^d(100) ) \approx (8, 6)$.

Also shown in  \Cref{f:Q2D_region+trajs} is the acceptance region obtained from an application of  Q-learning. 
Shown on the left-hand side is a sample path evolving within the decision boundary,   illustrating an instance where max-CUSUM* triggers a false alarm   avoided by the Q-learning policy.

\whamit{Performance evaluation.} For the results that follow we denote $\InfoState_k=(\InfoState_{k}^b,\InfoState_{k}^c)$ for \textbf{\modelCone} and $\InfoState_k=(\InfoState_{k}^*, \InfoState_{k}^d)$ for \textbf{\modelCtwo}.
\Cref{f:PerformanceQ2D} shows the cost obtained using Q-learning  as a function of $\kappa$.  Also shown for each case is  the cost using  CUSUM* for each of  the one-dimensional RRWs   defined by  the components of $\InfoState_k$.   For example,  $F^b=\surL+r^*$ with $\surL $ the LLR  defined with densities specified in \textbf{Model \modelAtwo}.

Most interesting are results obtained for \textbf{\modelCone}:
Q-learning  outperformed either CUSUM* cost for $\kappa\ge3$,  illustrating the value of using a more exotic SIS in applications with uncertain statistics.

\section{Conclusions} 
\label{s:conc}
	
	The theory and numerical results in this paper motivate many directions for future research:

\whamb
Further testing is needed beyond the independent model that is the focus of \Cref{s:Qqcd}---in particular,  in settings under which $\tchange$ and $\bfmX^0$ are not statistically independent as is the case in many applications.

\whamb
	In applications of interest to us there may be well understood behavior before a change  (which might represent a fault in a transmission line,  or a computer attack).    We cannot expect to have  a full understanding of post-change behavior.  The choice of SIS must be reconsidered in these settings,  perhaps based on techniques for universal hypothesis testing     (see \cite{huaunnmeyveesur11,huamey13} and the references therein). 
	
\whamb  
The MDP setting of \Cref{s:POMDPQCD} opens the door to stochastic gradient descent techniques to approximate the minimum of the  
 objective \eqref{e:ObjPOMDPQCD}.   Preliminary work in \cite{coomey24d} applied a version of the actor-critic method in application to QCD, but found the approach was not practical due to extraordinarily high variance of the gradient samples.   Variance reduction may be achievable through alternative approaches such as recent refinements of SPSA/extremum seeking control \cite{laumey25b}.  \notes{spa03 should be referenced in   dissertation}

\notes{the actor critic method is understood as a gradient-free method!  }
	
\clearpage

\bibliographystyle{abbrv}
\bibliography{strings,markov,q,extras}

\appendix

\section{Appendices}
\label{s:appndx}

\subsection{Convergence of Q-learning}
\label{s:SAthy}
	
	The proofs of 
\Cref{t:Q} and 
\Cref{t:Qgamma1}  are provided at the end of this subsection.   First we recall some general stochastic approximation (SA) theory.

\wham{Summary of general theory}   
	
The main results of \Cref{s:POMDP}
	are based in part on recent SA theory from  \cite{borchedevkonmey25}.    
	The standard $d$-dimensional SA recursion is expressed
\begin{equation}
	\theta_{n+1} =  \theta_n + \alpha_{n+1}  f(\theta_n,  \qsaprobe_{n+1})  \,,\quad  n\ge 0
	\label{e:SA}
\end{equation}
	In which $\theta_0\in\Re^d$ is given,  and $\{ \qsaprobe_{n+1} \}$ is a stochastic process---assumed here to be a Markov chain on  a general topological state space $\bigstate$, with Borel sigma algebra denoted $\bz$.   The  step-size sequence satisfies (A1) of  \cite{borchedevkonmey25}, which is denoted (SA1) here:
	
\wham{(SA1)}   The  step-size sequence is 
	non-negative,  satisfying   
\begin{equation} 
\begin{aligned}
			&\sum_{n=1}^{\infty} \alpha_n = \infty
\qquad \sum_{n=1}^{\infty} \alpha_n^2 < \infty
\\
			&
\lim_{n \to \infty} \big[ \alpha_{n+1}^{-1} - \alpha_n^{-1} \big]\ \ \textit{exists and is finite. } %\diffalpha = 
\end{aligned}
	\label{e:stepsize-cond}
\end{equation}
	
	Condition (SA1)  holds with  $\alpha_n = \min(\alpha_0,n^{-\rho})$ provided $\alpha_0>0$ and   $\half <\rho\le1$.

\wham{(SA2)}
	There exists a measurable function $\lips :\bigstate \to \posRe$ such that for each $z\in\bigstate$, 
\begin{equation}
	\label{e:lip-f}
\begin{aligned}
\| f(0, z)\| &\leq \lips(z), \\
\|f(\theta,z) - f(\theta', z)\| &\leq \lips(z)\|\theta - \theta'\| \,, \qquad \theta, \theta'\in \Re^d
\end{aligned}
\end{equation}

	The \textit{scaled vector field} is defined as $\displaystyle \barfinf(\theta)\eqdef \lim_{r\to\infty} \frac{1}{r}  \barf(r\theta) $, 
	and the   ODE@$\infty$ is then denoted $\displaystyle 
\ddt \odestate_t = \barf_{\infty}(\odestate_t)$.
	
\wham{(SA3)}   The mean flow is   exponentially asymptotically stable  (EAS) to some $\theta^*\in\Re^d$.   Moreover, the 
	limit defining the scaled vector field  exists for each $\theta\in\Re^d$ to define   a Lipschitz function $\barfinf\colon\Re^d\to\Re^d$.  
	The   ODE@$\infty$ is also EAS.

\wham{(SA4)} 
	The following drift condition holds:
	$$
\left. 
\mbox{\parbox{.85\hsize}{\raggedright
			For functions $V\colon\bigstate\to\Re_+$,  $ W\colon\bigstate\to [1, \infty)$  and $b>0$:\,
\[
\Expect\bigl[  \exp\bigr(  V(\qsaprobe_{k+1})      \bigr) \mid \qsaprobe_k=z \bigr]  
\le  \exp\bigr(  V(z)  - W(z) +  b   \bigl) \,, \qquad z\in\bigstate.
\]
	}}
\right\}
\eqno{\hbox{\bf (DV3)}}
	$$
	Moreover,  
	$\lips = o(W)$, meaning 
	$\displaystyle
\lim_{r\to\infty}  \sup_{z\in\bigstate}  \frac{ \lips(z) }{\max\{r,W(z)\}} =0 $, 
	and for  each $r>0$,
\begin{subequations}%
\begin{align}
\!   \!   \!   \!   \!
			S_W(r)  &:= \{ z :  W(z)\le r \}  \quad    \text{is either small or empty. }
\label{e:Wunbdd}
\\
			b_V(r) & \eqdef \sup\{ V(z) :  z\in S_W(r) \}  <\infty.
\label{e:V-bounded-level}
\end{align}
	\label{e:VWbdds}
\end{subequations}%

The following is obtained on combining Theorems~1 and 4 of \cite{borchedevkonmey25}.

\begin{proposition}[Convergence of SA]
\label[proposition]{t:SA}
If (SA1)--(SA4)  hold then the SA recursion 	\eqref{e:SA} is convergent with probability one from each initial condition $(\theta_0,\Phi_0)$.    

  Suppose that the step-size is
$ \alpha_{n + 1}  =1/ (n+1)^\rho$    with $\half < \rho < 1$.  Then
for the  estimates obtained via  averaging $\{\thetaPR_n \} $
as defined in  
\eqref{e:thetaPR}  we have,  as $n\to\infty$,
\[
\begin{aligned}
\zPR_n 
  & \distarrow  N(0, \SigmaPR )
\\
  n \Expect [ \tilthetaPR_n (\tilthetaPR_n)^\transpose  ] 
		& \longrightarrow	  \SigmaPR  
		\end{aligned}
\]
where
$\zPR_n\eqdef  \sqrt{n} \tilthetaPR_n$,
with $\tilthetaPR_n = \thetaPR_n -\theta^*$.
\qed
\end{proposition}

\wham{Implications to QCD theory}

	The Markov chain of interest in the proof of  \Cref{t:Q}  is the triple,  
\begin{equation}
	\qsaprobe_{k+1} = [U_k, \Phi_k,  \Phi_{k+1}  ]
	\label{e:Markov:tQ}
\end{equation}
	Under the assumption that $\bfmU = \{U_k \}$ is i.i.d.\ and independent of $\bfPhi$ it follows that $\bfqsaprobe$ is also Markovian.    
	Moreover, the recursion  \eqref{e:POMDPQ}  can be expressed in the form \eqref{e:SA} using \eqref{e:Markov:tQ}, in which 
	the function $f$ is described as follows:   Denoting $z=(u,x,x^+) \in \bigstate$ a possible value of $\qsaprobe_{k+1}$, and writing
	$s = g(x)$, $s^+ = g(x^+)$,
	\begin{equation}
		f(\theta,z) = \psi(s,u)
		\big[ - \theta^\transpose \psi(s,u) + c(x,u)
		+ \disc (1-u) \min_i \{ \theta^\transpose \psi(s^+,i) \} \big]
		\label{e:tQf}
	\end{equation}
	in which $\psi(s,u) = u \psi^1(s) + (1-u) \psi^0(s)$ for all $s\in\sstate$, $u\in\ustate$.

\Cref{t:QLip} explains the assumption  $G = o(W)$ imposed  in  \Cref{t:Q}.   Its proof 
	follows from \eqref{e:tQf} and
	elementary inequalities such as   $\|  \psi(s,u) \|  c(x,u)  \le \half [\|  \psi(s,u) \|^2 +  c(x,u) ^2 ]$.
	Recall  that $\InfoState_n = g(\Phi_n)$.

\begin{lemma}
	\label[lemma]{t:QLip}
		Under the assumptions of \Cref{t:Q}  the Lipschitz bound (SA2) holds with  
\begin{equation}
			\lips(z) =   K  \| \ell_\psi(x) \|^2    + c^2(x,u) + \| \ell_\psi(x^+) \|^2   
\label{e:LQCD}
\end{equation}
		in which   $\ell_\psi(x) = \|\psi^1(g(x))\| + \|\psi^0(g(x))\|$ for $x\in\state$, and    $K\ge 1$ is a sufficiently large constant.    
\end{lemma}

A greater challenge in the proof of \Cref{t:Q} is establishing (DV3) for $\{\qsaprobe_k \}$.
The following is required to obtain \eqref{e:VWbdds}.

\begin{lemma}
	\label[lemma]{t:QDV3a}
		Under the assumptions of \Cref{t:Q}  the  set $S_r$ below is either small or empty for any $r>0$:
	\[
		S_{W_1}(r) 
		= \{ z =  (u,x,x^+) \in \bigstate   :   W(x^+)  \le r \}
	\]
\end{lemma}
	
\wham{Proof}
	The Markov chains $\bfPhi$ and $\bfqsaprobe$ are each ``$\psi$-irreducible and aperiodic'' for their respective  probability measures $\psi$  (see \cite{MT} for definitions).
	Aperiodicity of  $\bfPhi$   implies that the set $S_W(r)$ is small if and only if there is a probability measure $\nu$ on $\bx$ and $\epsy>0$ such that
\begin{equation}
		K(x,A) \eqdef \sum_{k=1}^\infty 2^{-k} P^k(x,A)  \ge \epsy \nu(A)\,,\quad x\in S_W(r) \,, \  A\in\bx
	\label{e:SmallForPhi}
\end{equation}
	See \cite[Theorem 5.5.7]{MT}.
	The proof of the lemma requires that we establish an analogous bound for $\bfqsaprobe$.

	It is easy to show that the bivariate chain $(\Phi_k, \Phi_{k+1})$ admits an analogous bound with   
	probability measure $\nu_b$ defined for product sets via
\[
\nu_b(A\times B) = \int_{x\in A} \nu(x) P(x,B) \,,\quad A,B\in\bx
\]
	And then by independence of the input and $\bfPhi$ we have a candidate ``small'' probability measure for $\bfqsaprobe$:
\[
\nu_1(\{u\}\times A\times B) = \half  \nu_b(A\times B) \,,\quad A,B\in\bx\,, \, u\in\ustate
\]
	defined so that $\sum_u  \nu_1(\{u\}\times A\times B) =    \nu_b(A\times B) $.   If \eqref{e:SmallForPhi} holds then the following is immediate:
	For any $A_1 = \{u\}\times A\times B \in\bz$ and $z=(u,x^-,x) \in \bigstate$,  provided $W(x)\le r$,
\[
\begin{aligned}
		K_1(z, A_1)  &\eqdef \sum_{k=1}^\infty 2^{-k }  \Prob\{\qsaprobe_k\in A_1  \mid  \qsaprobe_0 = z\}
	\\
		&
		= \half  \sum_{k=1}^\infty 2^{-k }  \Prob\{ \Phi_k\in A\,, \Phi_{k+1} \in B  \mid    \Phi_0=x \}
	\ge  \nu_1(\{u\}\times A\times B)  
\end{aligned}
\]
\qed

\begin{lemma}
	\label[lemma]{t:QDV3b}
		Assumption (SA4) holds 
		under the assumptions of \Cref{t:Q} in the form 
	\[
	\Expect\bigl[  \exp\bigr(  V_1(\qsaprobe_{k+1})      \bigr) \mid \qsaprobe_k=z \bigr]  
	\le  \exp\bigr(  V_1(z)  - W_1(z) +  b_1   \bigl) \,, \qquad z\in\bigstate.
	\]
		in which  $V_1(\qsaprobe_{k+1}) = V(\Phi_{k+1}) + \half W(\Phi_{k}) $,
		$W_1(\qsaprobe_{k+1}) = 1+ \half  [ W(\Phi_{k+1})+ W(\Phi_{k})]$, and $b_1 = b+1$.  
		
		Moreover, 
		$V_1$,  $W_1$ satisfy   \eqref{e:VWbdds}, and 
		$\lips = o(W_1)$ with $\lips$ defined in  \eqref{e:LQCD}.
\end{lemma}
	
\wham{Proof}
	We   have $\lips = o(W_1)$ by the definition \eqref{e:LQCD}, and since the assumptions of the proposition imply that  
	$\| \ell_\psi \|^2 = o(W)$ and $  \ell_c^2 = o(W)$,  with $ \ell_c(x) = \max_i c(x,i)$.
	
	Next we establish (DV3).
	With   $z =  (u,x^-,x) \in \bigstate$ we have by the definitions,
\[
\begin{aligned}
	\Expect\bigl[  \exp\bigr(  V_1(\qsaprobe_{k+1})      \bigr) \mid \qsaprobe_k= z \bigr]  
		&\le  \exp\bigr(  V(x)    - \half W(x) +  b   \bigl)
	\\
		&= \exp\bigr(  V_1(z)    - \half [ W(x^-) + W(x) ]  +  b   \bigl) 
	\\
		&= \exp\bigr(  V_1(z)    - W_1(z)  +  b_1   \bigl) 
\end{aligned}
\]
	where the first equality uses $ V_1(z) = V(x) + \half W(x^-) $, and the final equation 
	uses $W_1(z) = 1+ \half  [ W(x^-) + W(x)]$ and 
	$b_1=b+1$.  
\qed
	
\wham{Proof of \Cref{t:Q}}
To apply \Cref{t:SA}, the verification of (SA1), (SA2), and (SA4) is provided in \Cref{t:QLip,t:QDV3a,t:QDV3b}.
It remains to verify (SA3), which involves two applications of \Cref{t:QstableMeanFlow}.   
First, the mean flow is EAS due to this lemma.
 
The vector field for the  ODE@$\infty$  is  $\barfinf(\theta) = \Expect_\piPhi[ f_\infty(\theta, \qsaprobe_k) ]$, where the expectation is in steady-state, 
\[
f_\infty(\theta, z) = \psi(s,u)\big[ -\theta^\transpose \psi(s,u) + \disc (1-u) \min_i\{ \theta^\transpose \psi(s',i) \} \big]
\]
where $z=(u,s,s')$.  Consequently, the vector field $\barfinf $  is precisely \eqref{e:barfQ} with $b=0$,  which is  \eqref{e:pBe} in the special case $c(x,u)\equiv 0$.     
Applying \Cref{t:QstableMeanFlow} we conclude that the ODE@$\infty$ is EAS as required.
\qed

\wham{Proof of \Cref{t:Qgamma1}}
	
	The proof is a minor modification of the proof of  \Cref{t:Qgamma1a}.
	
	Consider the positive kernel $P_\RegenState$ defined by $P_\RegenState (x,A) = \ind_\RegenState (x)   P(x,A)$  for $x\in\state$ and measurable $A\in\state$.    For $g\in L_2(\piPhi)$ we have, as in the proof of \Cref{t:Qgamma1a},   
\[
\| P_\RegenState g \|_\piPhi^2 = \| g \|_\piPhi^2  -\bigl( \sigma^g_{k+1\mid k}  +
\| \ind_{\RegenState^c} Pg \|_\piPhi^2  \bigr) 
\]
where $\sigma^g_{k+1\mid k}  $ is defined below  \eqref{e:Contract_a}.
Recall that  $ \| g \|_\piPhi^2 = \theta^\transpose R \theta$ if  $g= \theta^\transpose \psi$, giving
\[
\| P_\RegenState g \|_\piPhi^2 =      \theta^\transpose  \big[ R    -    \Sigma_{k+1\mid k}       -    M^\RegenState \big] \theta  
\]
	Under the assumptions of the proposition there is $\varrho>0$ such that 
\[
\| P_\RegenState g \|_\piPhi^2 
\le  \varrho \theta^\transpose R \theta = \varrho  \| g \|_\piPhi^2 
\]
	Hence $P_\RegenState$ is s a strict contraction on $\clH =\{  \theta^\transpose \psi : \theta\in\Re^d \}$,  and the proof follows as argued in 
\Cref{t:Qgamma1a}.
\qed
	
	%\clearpage

\subsection{Experiment details}
\label{s:appndx_exp}

\wham{Performance benchmarks}
For each model, Monte-Carlo simulation was used to estimate $\MDE(\thresh)$ and
$\MDD(\thresh)$ for CUSUM* and for Shiryaev's test
$p_n = \Prob\{\tchange \le n \mid \Obs_0^n\}$.  For CUSUM* we ran
$N = 2\times 10^6$ sample paths for a uniform grid
$0 < \thresh \le 20$ with $T=10^3$ thresholds; for Shiryaev's test we used a grid
$0 < \thresh \le 1$.  For each $\thresh$ the quantities $\MDE(\thresh)$ and
$\MDD(\thresh)$ were obtained as sample averages over the $N$ runs, yielding a
lookup table $\{\MDE(\thresh_t),\MDD(\thresh_t)\}_{t=1}^T$.  From these estimates
we evaluate the CUSUM* cost $\clJ(\thresh;\kappa)$ for $2\le \kappa\le 100$ and
approximate
\[
\thresh^*(\kappa) \eqdef \arg\min_{\thresh}\clJ(\thresh;\kappa),
\qquad
J^*(\kappa) \eqdef \clJ(\thresh^*(\kappa);\kappa),
\]
by minimizing over the grid.  These precomputed values define the CUSUM*
curves and the function $\clJ(\thresh;\kappa)$ used in the cost histograms.

\smallskip
\noindent
\wham{Q-learning}
Additional details of the Q-learning experiments are as follows:
\begin{itemize}
	\item[\wham{1.)}] \textit{Parameter initialization.}  
	For each run $1 \le m \le M$, we choose $\theta_0$ uniformly at random in $[-50, 50]$.
	
	\item[\wham{2.)}] \textit{Resetting.}  
	Large transients are common in stochastic approximation implementations.    In these
	experiments, whenever $\Vert \theta_k \Vert_\infty > 5\times 10^3$ we reset
	by resampling $\theta_{k+1}$ uniformly in $[-50, 50]$, independently of the
	past.
	
	\item[\wham{3.)}] \textit{Training horizon.}  
	All experiments to obtain $Q^{\theta^*}$ use $N = 2\times 10^5$ regeneration times,
	with a fixed constant $\upeta = 30$ used to define the regeneration set
	$\RegenState$ and the corresponding sequence $\{\Phi_k\}$.
	
	\item[\wham{4.)}] \textit{Basis parameters.}  
	Values $\{ \mu_i,\sigma_i \}_{i=1}^K$ in \eqref{e:basis2} were obtained by simulating $N = 2\times 10^4$ sample paths for each SIS. For each $1\le n\le N$ a common random seed was used in the initialization of all SIS trajectories;   and, as always, $\InfoState_k= 0$.
	
	The value $b=0.4$ was kept fixed in all Q-learning experiments. \textbf{QCD Models 1} and \textbf{2} used $K=20$, whereas \textbf{QCD Model 3} used $K=40$  to accommodate multidimensional features.

\end{itemize}
	
\end{document}